\documentclass[review,10pt]{elsarticle}
\usepackage[utf8]{inputenc}
\usepackage{microtype}
\usepackage[margin=2.5cm]{geometry}
\usepackage{setspace}
\usepackage[none]{hyphenat}
\usepackage[pagewise]{lineno}
\usepackage[english]{babel}
\usepackage{color}
\usepackage{booktabs}
\usepackage{longtable}
\usepackage[most]{tcolorbox}
\usepackage{ragged2e}
\usepackage[overload]{empheq}
\usepackage[justification=justified]{caption}
\usepackage{enumitem}

\usepackage{amsmath,amsfonts,amssymb,mathtools}
\usepackage{cases}
\usepackage{multirow}
\usepackage{graphicx,float}
\usepackage{algpseudocode}
\usepackage{subcaption}
\usepackage{eurosym}
\usepackage{amstext} 
\usepackage{arydshln}
\usepackage{flushend}
\usepackage{threeparttable}
\usepackage{url}
\usepackage{colortbl}
\usepackage{array}
\usepackage{multicol}
\usepackage{tabularx}
\usepackage{amsmath}
\usepackage{tikz}
\usepackage{optidef}
\definecolor{myblue}{RGB}{0,120,215}
\usetikzlibrary{positioning, calc, arrows.meta}

\usepackage[]{hyperref}
\hypersetup{colorlinks, linkcolor=blue, citecolor=blue, urlcolor=blue}

\usepackage[ruled,linesnumbered]{algorithm2e}
\DeclareRobustCommand{\officialeuro}{%
  \ifmmode\expandafter\text\fi
  {\fontencoding{U}\fontfamily{eurosym}\selectfont e}}
\usepackage{framed} 
\immediate\write18{%
makeindex Main_File.nlo -s nomencl.ist -o Main_File.nls 
}
\usepackage{multicol} 
\usepackage{nomencl}[noprefix] 
\makenomenclature
\renewcommand*\nompreamble{\begin{multicols}{2}}
\renewcommand*\nompostamble{\end{multicols}}

\def\tsc#1{\csdef{#1}{\textsc{\lowercase{#1}}\xspace}}
\tsc{WGM}
\tsc{QE}
\tsc{EP}
\tsc{PMS}
\tsc{BEC}
\tsc{DE}

\begin{document}

    \justifying
    \begin{frontmatter}   
        \title{DSO-Led Bilevel Optimization Framework for TSO–DSO Coordination across Active Distribution Networks}



        \author[1,2]{Fernando García Muñoz\texorpdfstring{\corref{mycorrespondingauthor}}{}} \cortext[mycorrespondingauthor]{Corresponding author}{\ead{fernando.garciam@usach.cl}}{} 

        \author[1]{Martín Venegas Escalona}

        \address[1]{University of Santiago of Chile (USACH), Faculty of Engineering, Industrial Engineering Department, Chile}

        \address[2]{University of Santiago of Chile (USACH), Faculty of Engineering, Program for the Development of Sustainable Production Systems (PDSPS), Chile}

    \begin{abstract}
         This work presents a bilevel coordination model that captures the hierarchical interaction between the transmission and distribution layers under a Distribution System Operator(DSO)-led configuration. In this scheme, multiple DSOs independently optimize the operation of their active distribution networks (ADNs), including photovoltaic (PV) generation, battery energy storage systems (BESS), and peer-to-peer (P2P) energy exchanges both within and across ADNs through the Transmission Network (TN), before the Transmission System Operator (TSO) performs the global coordination. The proposed formulation combines the Second-Order Cone relaxation of the DistFlow model to represent the distribution networks (DNs) with the classical DC optimal power flow (OPF) model for the transmission layer. The DSO-first decision sequence enables the reformulation of the bi-level problem into an equivalent single-level optimization model using the Karush–Kuhn–Tucker (KKT) conditions, resulting in a Mixed-Integer Second-Order Cone Programming (MISOCP) formulation that captures both the discrete and convex characteristics of the problem, while preserving the binary variables associated with DER and P2P operation, which would otherwise need to be relaxed in traditional TSO-led approaches. The model is tested on a hybrid system composed of the IEEE 30-bus transmission network and five IEEE 33-bus DNs. Results show that the DSO-led coordination leads to a more efficient use of BESS, improves local self-consumption, and reduces imports from the TN compared to the conventional top-down scheme. Furthermore, computational results from the case study reveal that the model exhibits near-linear or quadratic growth in problem size as the number of ADNs increases, suggesting its applicability to large-scale multi-ADN configurations. 
    \end{abstract}
    \begin{keyword}
    \emph{Bilevel optimization; Active distribution networks; Peer-to-peer energy trading}.
    \end{keyword}
\end{frontmatter}
\section{Introduction}
The global push for decarbonization, driven by the urgent need to reduce CO\textsubscript{2} emissions and limit global warming, has accelerated the deployment of renewable energy sources~\cite{PAPADIS2020118025, irena2024a}. As distributed energy resources (DERs), particularly solar PV and battery energy storage systems (BESS), become increasingly cost-competitive, power systems are undergoing a structural shift: distribution networks (DNs) are evolving from passive endpoints into active participants in energy provision and flexibility services, thus giving rise to the concept of Active Distribution Networks (ADNs) \cite{9328796}. This transformation has sparked growing interest in the coordination between Transmission System Operators (TSOs) and Distribution System Operators (DSOs), resulting in a diverse array of coordination schemes. In this regard, most existing studies have focused on schemes in which the TSO makes decisions first, while the DSO responds subsequently~\cite{lind2024tso, 9941459}. This sequential structure preserves the traditional top-down logic of power system operation, assigning a central role to the TSO and limiting the DSO to a reactive or subordinate position.

In this context, bilevel optimization has emerged as a prominent modeling approach, as it naturally captures the sequential and hierarchical structure of TSO-DSO coordination schemes. By explicitly distinguishing between upper-level (leader) and lower-level (follower) decisions, bilevel models reflect the reality of multi-actor power systems where the outcome of one agent's optimization problem influences the feasible set and objective of another. Nevertheless, under the conventional TSO-leader/DSO-follower paradigm, the subproblem associated with the DSO can become increasingly complex, particularly in ADNs with higher levels of decentralization and operational heterogeneity. The growing integration of electric vehicles, distributed generation, peer-to-peer (P2P) trading schemes, and local flexibility markets introduces numerous discrete decisions and nonlinear constraints~\cite{lind2019transmission}. These features give rise to high-dimensional, mixed-integer subproblems that are computationally challenging to solve. As a result, scalability becomes a critical bottleneck for the practical deployment of bilevel models following this traditional structure.

Alternative coordination schemes have emerged in response to these limitations, assigning DSOs a more central and proactive role as leaders in the decision-making process, while the TSO assumes a supporting or reactive role as the follower. This inversion of roles not only aligns better with the operational reality of increasingly decentralized ADNs but also could bring computational advantages. From a bilevel optimization perspective, placing the DSO at the upper level allows the most complex decisions, typically related to DER coordination, discrete operations, and local market participation, to be handled at the first stage, while the TSO-level problem, often focused on continuous and aggregated decisions, remains more tractable. In this context, the present work adopts a DSO-centric bilevel coordination scheme and proposes a mathematical model in which DSOs act as leaders managing local P2P energy exchanges and submitting aggregated information to the TSO. The TSO, as the follower, optimizes its operations based on the aggregated behavior of multiple DSOs and facilitates inter-network P2P exchanges across distribution areas.

The remainder of the paper is organized as follows. Section \ref{sec:Section_2} provides a literature review of existing coordination schemes, their corresponding optimization models, and solution strategies. Section \ref{sec:Section_3} introduces the main assumptions and the theoretical framework of the proposed coordination scheme. Section \ref{sec:Section_4} describes the mathematical formulation of the optimization models, details the solution approach, and presents the single-level reformulation. Section \ref{sec:Section_5} presents the case study and discusses the computational results. Finally, Section \ref{sec:Section_6} summarizes the main conclusions and outlines future research directions.

\section{Literature Review}\label{sec:Section_2}
Various coordination schemes have been proposed to facilitate the interaction between TSOs and DSOs in increasingly complex and decentralized power systems. While the literature presents a diverse array of frameworks, recent studies have highlighted the value of using Stackelberg formulations, both single-leader and multi-leader, as a theoretical foundation to interpret centralized versus decentralized coordination paradigms~\cite{YING2024123803}. These perspectives emphasize the importance of aligning the coordination structure with system architecture and computational tractability. In particular, the classification proposed in~\cite{lind2024tso} distinguishes between hierarchical and distributed coordination based on decision authority and information exchange, while~\cite{en15197312} offers a complementary synthesis focused on DER integration and market design. A similar systematization is presented in~\cite{GIVISIEZ2020106659}, where coordination models are categorized according to market maturity and grid architecture, and in~\cite{PAPADIS2020118025}, which identifies structural, regulatory, and operational barriers associated with different TSO-DSO configurations. Based on these reviews, coordination schemes can be broadly categorized into three groups~\cite{lind2024tso, en15197312, GIVISIEZ2020106659}: (\emph{i}) TSO-led schemes, where the DSO plays a reactive role and submits limited aggregated information upward; (\emph{ii}) distributed coordination schemes, which aim to preserve decentralization through iterative information exchanges or auxiliary market signals, without a clearly dominant actor; and (\emph{iii}) DSO-centric schemes, in which DSOs act as leaders, managing local DERs and market mechanisms before interacting with the TSO.

Most of the existing literature has focused on the traditional TSO-led coordination scheme, as classified earlier. This centralized structure aligns naturally with bilevel optimization models, which capture the sequential and hierarchical nature of the coordination process. Several works have applied this framework across different contexts. For instance,~\cite{YUAN2017600} develops a hierarchical coordination mechanism for joint energy and reserve dispatch, using a bilevel model in which the DSO submits generalized offer functions to the TSO to manage uncertainty and preserve local information. In~\cite{SOARES2020100333}, the authors propose a bilevel stochastic AC-OPF formulation to evaluate the provision of reactive power by DSOs under high DERs penetration, focusing on minimizing expected system costs while ensuring compliance with voltage and power flow constraints. The work in~\cite{BERALDOBANDEIRA2024110818} introduces an approximate dynamic programming model to aggregate distributed flexibility from DSOs, capturing operational constraints through radial network models based on DistFlow equations. In~\cite{9503337}, coordination is addressed through a Lagrangian-based decomposition method in which the TSO interacts with multiple DSOs, dynamically linearizing the transmission network (TN) to handle power flow nonlinearities. The authors in~\cite{9835136} address the heterogeneity across DSOs by modeling asymmetric reciprocal effects through scenario-based flexibility regions, allowing the TSO to anticipate operational constraints and prevent infeasibilities arising from unequal DER distributions.

In addition to the previously discussed contributions, several further studies have also adopted the traditional hierarchical coordination scheme. In~\cite{LIU202327}, the authors develop a robust transmission planning model coupled with a stochastic reinforcement solution for the distribution level, illustrating a sequential interaction where the TSO optimizes long-term investment decisions while DSOs adapt operationally under uncertainty. A similar structure is found in~\cite{CHEN2022118319}, which models the optimal participation of ADNs in energy and reserve markets, leveraging a two-stage stochastic program with TSO-DSO coupling through aggregated DER offers. In~\cite{HAJATI2024110840}, a local flexibility market is proposed in which the TSO sets the flexibility prices and the DSOs respond by adjusting local operations, preserving a centralized control structure. The review in~\cite{SOTO2021116268} further highlights that, despite growing interest in P2P mechanisms, most implementations still operate under TSO-dominant paradigms, limiting the autonomy of DSOs. For instance, in the Swiss case study presented in~\cite{KALANTARNEYESTANAKI2024110747}, flexibility exchanges are coordinated by the TSO, which specifies active and reactive power needs to which DSOs respond based on local capabilities. Authors in~\cite{RODRIGUES2023101204} formulate a reactive power management model where the TSO issues global voltage support requirements, and DSOs optimize local reactive injections accordingly. Additional studies have addressed specific TSO-led applications such as unit commitment~\cite{9006872, nawaz2020tso}, flexibility management~\cite{10609023, 7990560, 9939101, 9941459}, frequency regulation~\cite{1664986}, voltage control~\cite{9543104}, and market integration of distributed storage~\cite{4956966}.

The above bilevel models typically place the TSO in the upper level, shaping the feasible operational region of the DSOs, who act as followers optimizing local objectives under the constraints imposed by transmission-level decisions. To solve such models, two main methodological approaches have been adopted in the literature. The first is a mathematical reformulation of the lower-level DSO problem through Karush–Kuhn–Tucker (KKT) conditions \cite{zemkoho2021theoretical} or primal-dual transformations, resulting in a single-level mathematical program with equilibrium constraints (MPEC) \cite{dempe2015bilevel}. This strategy is commonly used in contexts where the DSO-level problem is convex and continuous, enabling the use of commercial solvers~\cite{SOARES2020100333, ZHANG2024123073, CHEN2022118319, RODRIGUES2023101204}. The second strategy involves decomposition-based algorithms, such as Benders decomposition~\cite{YUAN2017600, LIU202327}, the Column-and-Constraint Generation (CCG) algorithm~\cite{9835136}, or Lagrangian-based methods~\cite{9503337, KALANTARNEYESTANAKI2024110747}, which preserve the hierarchical nature of the original formulation and allow for scalable resolution in large-scale or multi-DSO settings.

While the aforementioned bilevel formulations under TSO-led coordination offer a rigorous framework for capturing the hierarchical nature of system operations, their scalability becomes increasingly challenging as ADNs incorporate more complex elements. Specifically, modeling DERs such as BESS, electric vehicles (EVs), demand response programs, and P2P trading mechanisms introduces additional binary variables and intertemporal constraints in the lower-level problem. Among the reviewed literature, only a few studies explicitly report acceptable computational performance at scale~\cite{LIU202327,9503337,9835136}; however, even these rely on solving subproblems with integer variables via direct mixed integer linear programming (MILP) formulations, which may compromise tractability in large-scale systems. As the role of DSOs expands and operational decisions become more decentralized and discrete, traditional hierarchical schemes may not suffice without incurring a significant computational burden. This has motivated recent efforts to explore the other two alternative coordination schemes.

As part of this ongoing research effort, a subset of the literature has focused on distributed coordination schemes, where the decision-making process is decoupled across system operators, and interactions between TSOs and DSOs are handled through iterative exchanges or market-based coupling mechanisms. These approaches aim to preserve data privacy, enhance modularity, and improve scalability in the presence of multiple DSOs or DER owners. For example,~\cite{10202840} proposes a distributed operational planning framework based on the Alternating Direction Method of Multipliers (ADMM), which enables coordinated scheduling of shared resources between TSO and DSOs while maintaining information confidentiality. Similarly,~\cite{9911669} introduces a transactive energy market model in which the clearing process is decomposed into subproblems for the TSO, DSOs, and DER agents, each solved independently and coordinated via coupling variables under an ADMM structure. A complementary contribution is found in~\cite{LI2024122328}, which integrates electricity and carbon trading mechanisms into a distributed coordination framework involving TSOs, DSOs, and prosumers. Similarly,~\cite{HUANG2022108179} explores a multi-aggregator coordination scheme in which local flexibility providers interact with the system operator through a hierarchical but decentralized exchange process. These models offer alternatives to hierarchical optimization, particularly for large-scale systems where decentralized decision-making and local autonomy are operationally and computationally advantageous. However, it is worth noting that these distributed schemes do not rely on bilevel optimization as a formal modeling framework. Instead, they adopt a hierarchical architecture resolved through decomposed or staged subproblems, where coordination is achieved iteratively rather than through explicit leader-follower formulations. While this enhances scalability and privacy, it may overlook the anticipative structure and strategic interdependence captured by bilevel models, which provide a stronger theoretical foundation for representing hierarchical decision-making and optimizing leader-follower interactions.

Building on the previous discussion, the third classification category, where DSOs act as leaders, emerges as a viable alternative. This configuration offers three main advantages: first, it enables modeling and managing the complexity of increasingly discrete ADNs at the upper level; second, it retains compatibility with bilevel optimization frameworks, unlike fully distributed approaches; and third, it allows for a more efficient allocation of resources, since decisions are made at the level where greater asset granularity and operational detail are available. In this regard, recent studies have begun to explore this perspective. For instance,~\cite{ZHANG2024123073} proposes a bilevel optimization model where the DSO manages demand-side flexibility at the upper level, while the TSO adjusts its dispatch decisions in response. A similar hierarchical structure is found in~\cite{9914683}, where DER and storage investment decisions are made at the distribution level and then passed to the TSO for system-wide operational planning. In~\cite{MARQUES2023101055}, the DSO leads the organization of local P2P flexibility markets, subsequently coordinating with the TSO to validate and integrate the resulting exchanges. Finally,~\cite{MANSOURI2023121062} introduces an interval-based nested framework in which the DSO aggregates flexibility from buildings and EV fleets, with the TSO reacting to these aggregated resources during system-level optimization.

Among the few previously reviewed studies that adopt a hierarchical structure centered on the DSO, only two, \cite{ZHANG2024123073} and \cite{9914683}, explicitly formulate a bilevel optimization model, while the remaining two follow a sequential, multi-stage approach. None of these works incorporates P2P energy trading mechanisms among prosumers, thereby reducing the operational complexity faced by the DSO in managing centralized assets such as BESS, typically modeled with a limited number of binary variables. Moreover, all of them assume a single DSO, overlooking the interactions that naturally emerge across multiple DNs. This simplification not only limits the realism of the proposed frameworks, since a TSO typically coordinates with several DSOs, but also compromises their scalability when applied to larger systems. In this regard, this work continues the development of DSO-led coordination schemes by addressing part of the research gaps identified above through the following contributions: 

\begin{itemize}
    \item A bilevel optimization framework is proposed to represent the hierarchical coordination between multiple ADNs and the TN, where DSOs act as leaders and a common TSO acts as the follower. Each DSO optimizes the operation of local DERs and P2P trading within its ADN, while the TSO subsequently coordinates inter-ADN exchanges and manages transmission-level generation, including low-cost PV and high-cost conventional units, to ensure global feasibility.
    \item A tractable single-level MISOCP formulation is derived by embedding the TSO problem into the upper-level model through its KKT optimality conditions. Leveraging the convexity of the TSO’s DC-OPF problem, this reformulation preserves the hierarchical semantics of the bilevel structure while retaining the binary operational variables associated with DER and P2P decisions that are typically relaxed in conventional TSO-led approaches.
    \item The proposed framework enables the analysis of how a DSO-led coordination scheme could lead to a more efficient management of BESS resources at the distribution level compared to a traditional TSO-led configuration, capturing the associated impacts on local self-consumption and distributed flexibility. Furthermore, its ADN-centered structure allows scalable coordination, as the formulation exhibits a near-linear growth in computational complexity with the number of connected networks.
\end{itemize}


\section{Coordination scheme}\label{sec:Section_3}
This section introduces the proposed coordination scheme in which DSOs act as leaders and the TSO acts as the follower. It first describes the main operational assumptions underlying this DSO-led structure, outlining the hierarchical sequence of decisions and the interaction between transmission and distribution layers. Then, it presents the theoretical framework that formalizes this coordination as a bilevel optimization problem, providing the basis for the subsequent mathematical formulation of the model.

The proposed model represents a day-ahead coordination framework with hourly time steps, in which different DSOs independently operate multiple ADNs. Each ADN hosts local DERs, namely PV generation and BESS, enabling P2P energy trading among prosumers. P2P transactions are allowed within the same ADN and across different ADNs, via the transmission-level coordination of the TSO. In parallel, the TN includes utility-scale generation resources managed by the TSO, including low-cost PV generation units and higher-cost conventional generators. 

We assume the existence of a robust regulatory and technical infrastructure that supports secure P2P transactions, transparent exchange of information, and hierarchical coordination between DSOs and the TSO. All P2P transactions occur under a trading price known beforehand and lie between the purchase price from the TSO and the compensation rate for surplus injections into the grid. Within this structure, each DSO autonomously optimizes the operation of its local resources over a 24-hour horizon, aiming to maximize self-consumption, minimize operational costs, and coordinate P2P energy exchanges. After internal optimization, each DSO submits its net energy exchange profile, aggregating local demand, DER output, storage operation, and P2P trading, to the TSO. The TSO subsequently optimizes the allocation of centralized transmission-level resources in response to the aggregated exchanges offered by the DSOs. 

\begin{figure}[htbp]
    \centering
    \includegraphics[width=0.8\linewidth]{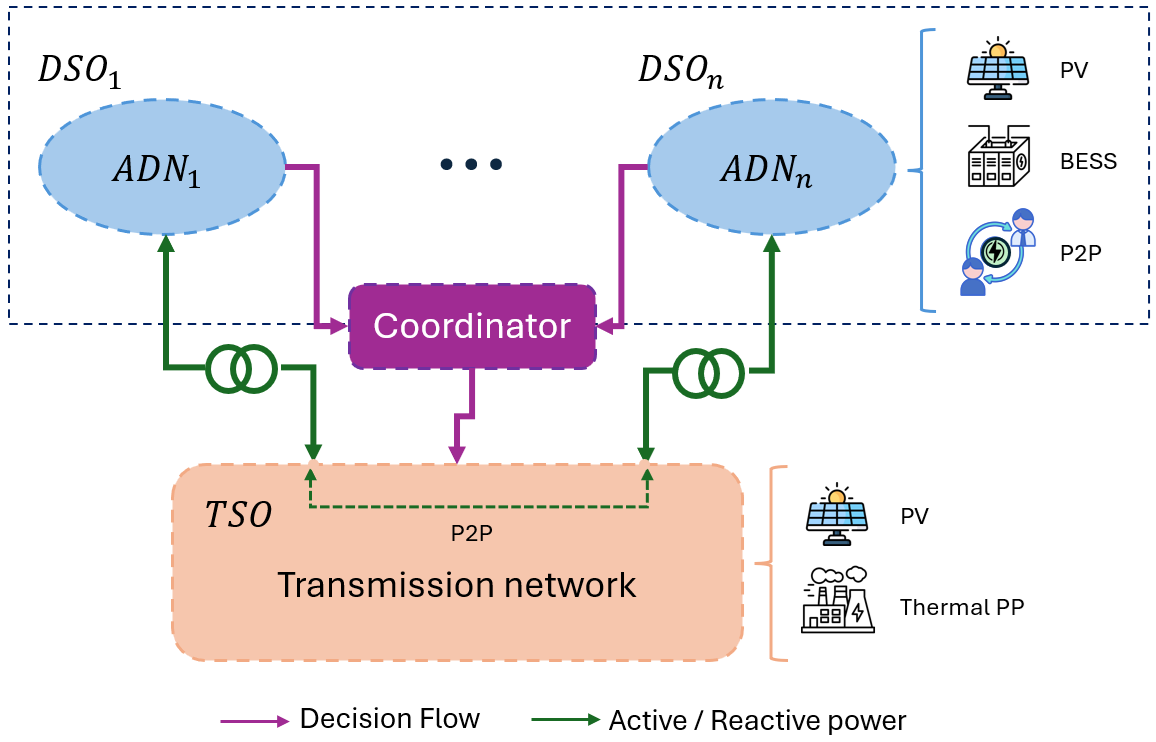}
    \caption{Schematic representation of the hierarchical connection between the TN and multiple ADNs}
    \label{fig:test_system_architecture}
\end{figure}

The availability of limited low-cost PV generation at the transmission level introduces an implicit interaction among DSOs, as they independently seek to minimize their operational costs while accessing this shared resource. In this regard, coordination is implicitly achieved through the centralized decisions of the TSO, which acts as a balancing agent. Thus, the DSOs are assumed to operate under symmetric regulatory conditions and do not engage in strategic behavior toward one another. In other words, while each DSO makes autonomous decisions based on local information and system parameters, it does not anticipate or react to the actions of other DSOs. This abstraction captures an operational paradigm in which DSOs follow standardized procedures and interact only indirectly through the central coordinator. Note that, since each DSO is effectively solving an economic dispatch problem, the generation cost of the TN is incorporated into its objective function as a cost-based price signal for accessing that shared resource, rather than as a market-clearing electricity price. Therefore, the proposed model, summarized in Figure~\ref{fig:test_system_architecture}, adopts an inverted coordination scheme in which DSOs act first by independently optimizing the operation of their local DERs, while the TSO subsequently optimizes the dispatch of transmission-level resources based on the aggregated information received from all DSOs.


Building upon these operational assumptions, the coordination problem between DSOs and the TSO can be represented through a bilevel optimization framework, where the upper level corresponds to the DSOs' aggregated decision-making process and the lower level represents the TSO’s operational response. Although multiple ADNs are connected to the TN, their actions are assumed to be centrally coordinated and non-strategic, reflecting the existence of common regulatory frameworks and standardized operational procedures. Under this assumption, the DSOs’ collective behavior can be represented as a single decision-making entity, which acts as the leader in a hierarchical structure, while the TSO acts as the follower optimizing its dispatch in response to the aggregated exchanges from the distribution level. This setting corresponds to a Single-Leader Single-Follower (SLSF) Stackelberg structure~\cite{zemkohoo2020bilevel}, which captures the hierarchical nature of the coordination while preserving a separation between local and system objectives.

Formally, the upper-level problem aggregates the decisions of all DSOs regarding local operation, DER coordination, and boundary energy exchanges, anticipating the optimal reaction of the TSO represented by the lower-level problem. The general structure can be expressed as:
\begin{equation}\label{SLSF_Form}
\begin{aligned}
    \min_{x, y} \quad & F(x, y) = \sum_{\nu=1}^{N} F(x, y) \\
    \text{s.t.} \quad & x \in X \\
    & y \in \arg\min_{y} \{G(y; x) \;\mid\; y \in Y(x) \}
\end{aligned}
\end{equation}

Where $x$ denotes the collective decisions at the distribution level and $y$ the optimal response of the TN. This Stackelberg representation is consistent with the coordination assumptions introduced above: the DSOs decide first, optimizing the operation of their ADNs under local and regulatory constraints, and the TSO subsequently adjusts the operation of transmission resources based on these aggregated boundary exchanges. In contrast to the top-down configuration, where the TSO optimizes based on statistically aggregated representations of distribution-level behavior, the Stackelberg formulation adopted here allows the TSO to respond to aggregated outcomes derived from individual DSO decisions. This distinction is essential: while both layers exchange aggregated information, in the DSO-first sequence, such aggregation reflects optimized local behaviors rather than assumed averages, leading to a more coherent and operationally meaningful coordination outcome. Consequently, the coordination scheme addressed in this work, and modeled through a bilevel Stackelberg framework, not only captures the physical coupling between transmission and DNs but also manages the granularity of individual DER decisions at the distribution level.

\section{Optimization models}\label{sec:Section_4}
This section presents the mathematical formulations used to represent the operation of the TN and the ADNs, as part of the hierarchical TSO-DSO coordination framework. All symbols, variables, and parameters are defined in the text when first introduced. A complete list of the notation employed in the models is provided in Table~\ref{Nomenclature_Table} for quick reference.

\subsection{Transmission Network}\label{Subsec:TN_C}
TNs operate under AC conditions, where power flows and voltages are governed by nonlinear equations involving active and reactive power. However, for operational planning and dispatch coordination, the TN can be approximated under assumptions valid in transmission systems such as voltage magnitudes close to 1~p.u., small phase angle differences, and negligible reactive power flows~\cite{4956966,cain2012history}. Under these conditions, the AC power flow equations can be linearized without significantly compromising accuracy, resulting in a tractable and computationally efficient model~\cite{wood2013power}, which is known as the DC Optimal Power Flow (DC-OPF), and is adopted in this work to represent the transmission conditions.

In the DC-OPF framework, the TN is represented as an undirected graph \( G = (\Omega_T, \mathcal{L}_T) \), where \( \Omega_T \) denotes the set of transmission buses and \( \mathcal{L}_T \) the set of transmission lines. Each line \( (i,j) \in \mathcal{L}_T \) connects two buses and supports bidirectional active power flows. A subset of buses \( \Omega^T_B \subset \Omega_T \) is designated as boundary nodes, serving as interconnection points between the TN and the DNs, each of which is independently operated by a DSO. The TSO minimizes in~\eqref{objfuncTSO} the total operational generation cost by coordinating two distinct types of generation technologies: conventional thermal generation and PV solar generation. Conventional units are modeled with a standard quadratic cost function that captures fuel consumption using coefficients \( Ca^T_i \), \( Cb^T_i \), and \( Cc^T_i \) for each TN bus \( i \). In contrast, PV generation considers lower operational cost, which is reflected through a linear term \( \pi^T_i pv^T_{i,t} \), where \( \pi^T_i \) represents a marginal cost associated with the solar generation at bus \( i \). The complete objective function is given by:
\begin{flalign} 
\label{objfuncTSO}
&\text{min} \enspace z = \displaystyle\sum_{t \in \mathcal{T}} \displaystyle\sum_{i\in \Omega_T} (Ca^T_{i} (pg^T_{i,t})^2 + Cb^T_{i} (pg^T_{i,t}) + Cc^T_{i}) + \displaystyle\sum_{t \in \mathcal{T}} \displaystyle\sum_{i\in \Omega_T}\pi^T_{i} (pv^T_{i,t})
\end{flalign}

The TN operation is subject to the constraints in~\eqref{TSO}, where power balance at each bus is enforced through~\eqref{TSOa} depending on the type of bus. For internal buses \( i \in \Omega^T \setminus \Omega^T_B \), the net injection equals the difference between conventional generation \( pg^T_{i,t} \) and PV generation \( pv^T_{i,t} \), minus the local demand \( PL^T_{i,t} \). At boundary buses \( i \in \Omega^T_B \), the balance accounts for energy exchanges with the connected ADNs, represented by the variables \( pk^{T,bg}_{i,t} \), \( pk^{T,sgc}_{i,t} \), and \( pk^{T,sge}_{i,t} \), which denote the power purchased from ADNs, and the power sold to ADNs as cheap and expensive energy blocks, respectively. In this regard, the amount of power injected into ADNs is constrained by~\eqref{TSOb} and~\eqref{TSOc}, which limit the maximum sales from PV and conventional generation, respectively. Additionally,~\eqref{TSOd} imposes an upper bound \( K^{T,bg}_i \) on the power that can be purchased from ADNs at each boundary bus.

Constraint~\eqref{TSOe} models the active power flow \( p^T_{i,j,t} \) on each line \( (i,j) \), assuming it is proportional to the difference in voltage angles \( \theta^T_{i,t} - \theta^T_{j,t} \) and inversely proportional to the line reactance \( X^T_{i,j} \). The resulting flow is limited by the thermal capacity of the line, as imposed by~\eqref{TSOf}, with upper bound \( S^{T,\max}_{i,j} \). Constraint~\eqref{TSOg} restricts active power generation at each bus to the maximum available generation \( PG^{T,\max}_{i,t} \). Likewise, PV generation at each bus is limited by Constraint~\eqref{TSOh}, which depends on the maximum available PV generation \( PV^{T,\max}_{t} \) and the installed PV capacity ratio \( \gamma^{T,pv}_{i} \). Finally, voltage angle references are established at selected reference buses \( \Omega_{ref} \) by enforcing Constraint~\eqref{TSOi}, where voltage angles are fixed at zero. The resulting set of operational constraints is defined as follows:

\begin{subequations}\label{TSO}
\begin{align} 
    &\displaystyle\sum_{(i,j) \in \mathcal{L}_T} p^T_{i,j,t} = 
    \left \{
  \begin{aligned}
    &pg^T_{i,t} + pv^T_{i,t} - PL^T_{i,t} && \forall i \in \Omega_T\\
    &pk^{T,bg}_{i,t} - pk^{T,sgc}_{i,t} - pk^{T,sge}_{i,t} && \forall i \in \Omega^T_{B} 
  \end{aligned} \right. && \forall t \in \mathcal{T} \label{TSOa}\\
    &\displaystyle\sum_{i \in \Omega^T_{B}}pk^{T,sgc}_{i,t} \leq \displaystyle\sum_{i \in \Omega_T} pv^{T}_{i,t}  && \forall t \in \mathcal{T} \label{TSOb}\\
    &\displaystyle\sum_{i \in \Omega^T_{B}}pk^{T,sge}_{i,t} \leq \displaystyle\sum_{i \in \Omega_T} pg^{T}_{i,t} && \forall t \in \mathcal{T} \label{TSOc}\\
    & pk^{T,bg}_{i,t} \leq K^{T,bg}_{i} && \forall i \in \Omega^T_{B}, \forall t \in \mathcal{T} \label{TSOd}\\
    &p^T_{i,j,t} = \frac{\theta^T_{i,t} - \theta^T_{j,t}}{X^T_{i,j}} && \forall (i,j) \in \mathcal{L}_T, \, \forall t \in \mathcal{T}  \label{TSOe}\\
    &- S^{T,\max}_{i,j} \leq p^T_{i,j,t} \leq S^{T,\max}_{i,j} && \forall (i,j) \in \mathcal{L}_T, \, \forall t \in \mathcal{T}  \label{TSOf}\\
    &pg^T_{i,t} \leq PG^{T,\max}_{i,t} && \forall i \in \Omega_T, \, \forall t \in \mathcal{T} \label{TSOg}\\
    &pv^T_{i,t} \leq PV^{T,\max}_{t} \gamma^{T,pv}_{i} && \forall i \in \Omega_T, \, \forall t \in \mathcal{T} \label{TSOh}\\
    & \theta^T_{i,t} = 0 && \forall i \in \Omega_{ref}, \, \forall t \in \mathcal{T} \label{TSOi} 
\end{align}
\end{subequations}

\subsection{Distribution Network}\label{Subsec:DN_C}
Each ADN \( k \in \Omega^T_B \) is operated by a DSO and modeled using a second-order cone relaxation of the AC-OPF~\cite{1664986,9543104}. The DN is represented by a directed graph \( G = (\Omega_{D_k}, \mathcal{L}_{D_k}) \), where \( \Omega_{D_k} \) denotes the set of buses in distribution system \( k \), and \( \mathcal{L}_{D_k} \subseteq \Omega_{D_k} \times \Omega_{D_k} \) defines the set of distribution lines connecting adjacent nodes. This formulation relies on the branch flow model~\cite{7990560}, which accurately represents voltage magnitudes, power flows, and current relationships in radial topologies, while enabling a convex relaxation that ensures computational tractability. The ADN includes DERs, namely PV generation and BESS, installed at specific nodes \( i \in \Omega^D_{A_k} \), which also participate in local P2P energy trading within and across distribution systems. The DSO aims to minimize the total energy procurement cost by coordinating distributed resources and boundary exchanges with the TN. The terms \( pk^{D,bgc}_{i,t} \), \( pk^{D,bge}_{i,t} \), and \( pk^{D,sg}_{i,t} \) represent, respectively, the power bought from the TN at cheap and expensive blocks, and the power injected into the TN from DN \( k \) at boundary bus \( i \) during period \( t \). These exchanges are valued according to the generation cost of the TN, incorporated into the DSOs’ objective function as cost-based signals, $\lambda^{D,bgc}_t$, $\lambda^{D,bge}_t$, and $\lambda^{D,sg}_t$, rather than as market-clearing electricity tariffs. 
The resulting objective function~\eqref{objfuncDSO} for the DSO is given by:
\begin{flalign}
\label{objfuncDSO}
&\text{min} \enspace z = 
\displaystyle\sum_{t \in \mathcal{T}} \displaystyle\sum_{k\in \Omega^T_{B}} \displaystyle\sum_{i\in \Omega_{D_k}} (\lambda^{D,bgc}_{t} pk^{D,bgc}_{i,t} + \lambda^{D,bge}_{t}  pk^{D,bge}_{i,t} - \lambda^{D,sg}_{t} pk^{D,sg}_{i,t})
\end{flalign}

The physical operation of each ADN is governed by the set of constraints~\eqref{DSO}. Active power balance is enforced in~\eqref{DSOa}, differentiating between regular agents \( i \in \Omega^D_{A_k} \) and boundary buses \( i \in \Omega^D_{B_k} \). The net injection at each node considers active power flows \( p^D_{i,j,t} \) between buses \( i \) and \( j \), and ohmic losses proportional to the squared current \( \ell^D_{j,i,t} \) and line resistance \( R^D_{j,i} \). Reactive power balance is imposed in~\eqref{DSOb}, where flows \( q^D_{i,j,t} \) and losses \( X^D_{j,i} \ell^D_{j,i,t} \) are balanced with the reactive generation \( qg^D_{i,t} \) and the reactive demand \( QL^D_{i,t} \). Constraints~\eqref{DSOc} defines the nodal active power balance for regular agents, including PV generation \( pv^D_{i,t} \), active demand \( PL^D_{i,t} \), and battery operation via charging \( ch^{D,bt}_{i,t} \) and discharging \( ds^{D,bt}_{i,t} \). At boundary buses, the net power exchange is represented in~\eqref{DSOd} by the variables \( pk^{D,bgc}_{i,t} \), \( pk^{D,bge}_{i,t} \), and \( pk^{D,sg}_{i,t} \), which denote, respectively, the power imported from the TN using cheap and expensive blocks, and the power exported to the TN.

Voltage magnitudes \( v^D_{i,t} \) are related through the branch flow model in~\eqref{DSOe}, where voltage drop depends on line parameters \( R^D_{i,j} \), \( X^D_{i,j} \), and power flows. The convex relaxation is enforced in~\eqref{DSOf}, where the apparent power squared is upper bounded by the product of the sending-end voltage and squared current magnitude. Limits on reactive generation \( qg^D_{i,t} \) are imposed in~\eqref{DSOg}, bounded by \( QG^{D,\min}_i \) and \( QG^{D,\max}_i \). Voltage magnitudes are constrained in~\eqref{DSOh} within the operational range defined by \( V^{D,\min}_i \) and \( V^{D,\max}_i \). Lastly, the squared current \( \ell^D_{i,j,t} \) is upper bounded in~\eqref{DSOi} by the thermal limit \( I^{D,\max}_{i,j} \), and Constraints~\eqref{DSOp2pa}–\eqref{DSOp2pc} set upper bounds on the power that can be sold or purchased at boundary buses, distinguishing between cheap and expensive energy blocks.
\begin{subequations}\label{DSO}
\begin{align}
    &\displaystyle\sum_{(i,j) \in \mathcal{L}_{D_k}} p^D_{i,j,t} - \displaystyle\sum_{(j,i) \in \mathcal{L}_{D_k}} (p^D_{j,i,t} - R^D_{j,i}\ell^D_{j,i,t}) =  
    \left \{
      \begin{aligned}
        &\Delta p^D_{i,t} && \forall i \in \Omega^D_{A_k} \ \\
        &\Delta p^{D,B}_{i,t} && \forall i \in \Omega^D_{B_k} \
      \end{aligned} \right. && \forall t \in \mathcal{T} \label{DSOa}\\
    &\displaystyle\sum_{(i,j) \in \mathcal{L}_{D_k}} q^D_{i,j,t} - \displaystyle\sum_{(j,i) \in \mathcal{L}_{D_k}} (q^D_{j,i,t} - X^D_{j,i}\ell^D_{j,i,t}) =  qg^D_{i,t} - QL^D_{i,t} && \forall i \in \Omega_{D_k}, \, \forall t \in \mathcal{T} \label{DSOb}\\
    &\Delta p^D_{i,t} = pv^D_{i,t} - PL^D_{i,t} - ch^{D,bt}_{i,t} + ds^{D,bt}_{i,t} && \forall i \in \Omega^D_{A_k}, \, \forall t \in \mathcal{T} \label{DSOc}\\
    &\Delta p^{D,B}_{i,t} = pk^{D,bgc}_{i,t} + pk^{D,bge}_{i,t} - pk^{D,sg}_{i,t} && \forall i \in \Omega^D_{B_k}, \, \forall t \in \mathcal{T} \label{DSOd}\\
    &v^D_{j,t} = v^D_{i,t} -2 (R^D_{i,j} p^D_{i,j,t} + X^D_{i,j} q^D_{i,j,t}) + ((R^D_{i,j})^2 + (X^D_{i,j})^2)\ell^D_{j,i,t} && \forall (i,j) \in \mathcal{L}_{D_k}, \, \forall t \in \mathcal{T} \label{DSOe}\\
    &(p^D_{i,j,t})^2 + (q^D_{i,j,t})^2 \leq \ell^D_{j,i,t} v^D_{i,t} && \forall (i,j) \in \mathcal{L}_{D_k}, \, \forall t \in \mathcal{T}
     \label{DSOf}\\
    &QG^{D,\min}_{i} \leq qg^D_{i,t} \leq QG^{D,\max}_{i} && \forall i \in \Omega_{D_k}, \, \forall t \in \mathcal{T} \label{DSOg}\\
    &V^{D,\min}_{i} \leq v^D_{i,t} \leq V^{D,\max}_{i} && \forall i \in \Omega_{D_k}, \, \forall t \in \mathcal{T} \label{DSOh}\\
    &\ell^D_{i,j,t} \leq I^{D,\max}_{i,j} && \forall (i,j) \in \mathcal{L}_{D_k}, \, \forall t \in \mathcal{T} \label{DSOi}\\
    &pk^{D,sg}_{i,t} \leq K^{D,sg}_{i} && \forall i \in \Omega^D_{B_k}, \, \forall t \in \mathcal{T} \label{DSOp2pa}\\
    &pk^{D,bgc}_{i,t} \leq K^{D,bgc}_{i} && \forall i \in \Omega^D_{B_k}, \, \forall t \in \mathcal{T} \label{DSOp2pb}\\
    &pk^{D,bge}_{i,t} \leq K^{D,bge}_{i} && \forall i \in \Omega^D_{B_k}, \, \forall t \in \mathcal{T} \label{DSOp2pc}
\end{align}
\end{subequations}

PV power output at node \( i \in \Omega^D_{A_k} \) and time \( t \), denoted by \( pv^D_{i,t} \), is constrained in~\eqref{DSOpv} by the product of the irradiance-dependent availability \( PV^{D,\max}_{t} \) and the installed capacity ratio \( \gamma^{D,pv}_{i} \), which reflects the fraction of maximum installable capacity deployed at node \( i \). BESS are modeled through a set of operational constraints~\eqref{DSObt}, where the state of charge (SoC), denoted by \( soc^{D,bt}_{i,t} \), evolves according to~\eqref{DSObta}, based on the charging power \( ch^{D,bt}_{i,t} \), discharging power \( ds^{D,bt}_{i,t} \), efficiencies \( \varphi^{D,ch} \) and \( \varphi^{D,ds} \), and the previous SoC. Constraint~\eqref{DSObtb} bounds the SoC between minimum and maximum permissible levels \( SOC^{D,bt}_{\min} \) and \( SOC^{D,bt}_{\max} \), scaled by the installed BESS capacity ratio \( \gamma^{D,bt}_{i} \). Charging and discharging powers are limited in~\eqref{DSObtc}–\eqref{DSObtd} by the rated power \( PB^{D,bt} \) and the binary variable \( w^{D,bt}_{i,t} \), which indicates the battery's operational mode. The binary parameter \( \nu^{D,bt}_{i} \) identifies whether a battery is installed at node \( i \), ensuring that discharging is disabled when no storage is available. Constraints~\eqref{DSObte} enforce logical consistency by allowing battery operation only if it is installed.

\begin{align}
    &pv^D_{i,t} \leq PV^{D,\max}_{t} \gamma^{D,pv}_{i} && \forall i \in \Omega^D_{A_k}, \, \forall t \in \mathcal{T} \label{DSOpv}
\end{align}

\begin{subequations}\label{DSObt}
\begin{align}
    &soc^{D,bt}_{i,t+1} = soc^{D,bt}_{i,t} + (\varphi^{D,ch} ch^{D,bt}_{i,t} - \frac{1}{\varphi^{D,ds}} ds^{D,bt}_{i,t}) \Delta t && \forall i \in \Omega^D_{A_k}, \, \forall t \in \mathcal{T} \label{DSObta}\\
    &SOC^{D,bt}_{\min} \gamma^{D,bt}_{i} \leq soc^{D,bt}_{i,t} \leq SOC^{D,bt}_{\max} \gamma^{D,bt}_{i} && \forall i \in \Omega^D_{A_k}, \, \forall t \in \mathcal{T} \label{DSObtb}\\
    &ch^{D,bt}_{i,t} \leq  PB^{D,bt} (w^{D,bt}_{i,t}) && \forall i \in \Omega^D_{A_k}, \, \forall t \in \mathcal{T} \label{DSObtc}\\
    &ds^{D,bt}_{i,t} \leq  PB^{D,bt} (1 - w^{D,bt}_{i,t}) - PB^{D,bt} (1 - \nu^{D,bt}_{i}) && \forall i \in \Omega^D_{A_k}, \, \forall t \in \mathcal{T} \label{DSObtd}\\
    &w^{D,bt}_{i,t} \leq \nu^{D,bt}_{i} && \forall i \in \Omega^D_{A_k}, \, \forall t \in \mathcal{T} \label{DSObte}
\end{align}
\end{subequations}

The P2P trading model allows agents within the same DN to exchange energy locally or through the TN interface. Each agent’s net traded power \( \Delta p^D_{i,t} \) is defined in~\eqref{DSOp2pd} as the difference between energy supplied \( \Delta p^{D,+}_{i,t} \) and energy demanded \( \Delta p^{D,-}_{i,t} \). If an agent experiences an energy deficit (\( \Delta p^D_{i,t} < 0 \)), it can purchase energy either from the local market (\( p^{D,bm}_{i,t} \)) or from the TN (\( p^{D,bg}_{i,t} \)). Conversely, in case of an energy surplus (\( \Delta p^D_{i,t} > 0 \)), the agent can sell excess energy to local peers (\( p^{D,sm}_{i,t} \)) or inject it into the TN (\( p^{D,sg}_{i,t} \)), as modeled in~\eqref{DSOp2pe}–\eqref{DSOp2pf}. To prevent agents from simultaneously acting as both buyers and sellers at a given time \( t \), binary variables \( y^D_{i,t} \) are introduced, where \( y^D_{i,t} = 1 \) indicates that agent \( i \) acts as a seller at time \( t \). Constraints~\eqref{DSOp2pg}–\eqref{DSOp2ph} enforce this exclusivity by allowing either selling or buying actions, but not both, within the same time period. The constant \( \mathcal{M}^D \) is a sufficiently large positive scalar used to activate or deactivate the corresponding decision variables. The local market equilibrium conditions are represented through energy balance constraints~\eqref{DSOp2pi}–\eqref{DSOp2pk}, which guarantee that total local sales match total local purchases and that total boundary exchanges are consistent with the aggregation of individual trades.

\begin{subequations}\label{DSOp2p}
\begin{align}
    &\Delta p^D_{i,t} = \Delta p^{D,+}_{i,t} - \Delta p^{D,-}_{i,t} && \forall i \in \Omega^D_{A_k}, \, \forall t \in \mathcal{T} \label{DSOp2pd}\\
    &\Delta p^{D,+}_{i,t} = p^{D,sg}_{i,t} + p^{D,sm}_{i,t} && \forall i \in \Omega^D_{A_k}, \, \forall t \in \mathcal{T} \label{DSOp2pe}\\
    &\Delta p^{D,-}_{i,t} = p^{D,bg}_{i,t} + p^{D,bm}_{i,t} && \forall i \in \Omega^D_{A_k}, \, \forall t \in \mathcal{T} \label{DSOp2pf}\\
    &\Delta p^{D,+}_{i,t} \leq \mathcal{M}^D (y^D_{i,t}) && \forall i \in \Omega^D_{A_k}, \, \forall t \in \mathcal{T} \label{DSOp2pg}\\
    &\Delta p^{D,-}_{i,t} \leq \mathcal{M}^D (1-y^D_{i,t}) && \forall i \in \Omega^D_{A_k}, \, \forall t \in \mathcal{T} \label{DSOp2ph}\\
    &\displaystyle\sum_{i \in \Omega^D_{A_k}} p^{D,sm}_{i,t} = \displaystyle\sum_{i \in \Omega^D_{A_k}} p^{D,bm}_{i,t} && \forall t \in \mathcal{T} \label{DSOp2pi}   \\
    &\displaystyle\sum_{i \in \Omega^D_{A_k}} p^{D,sg}_{i,t} = \displaystyle\sum_{i \in \Omega^D_{B_k}} pk^{D,sg}_{i,t} && \forall t \in \mathcal{T}\label{DSOp2pj}\\
    &\displaystyle\sum_{i \in \Omega^D_{A_k}} p^{D,bg}_{i,t} = \displaystyle\sum_{i \in \Omega^D_{B_k}} pk^{D,bgc}_{i,t} + \displaystyle\sum_{i \in \Omega^D_{B_k}} pk^{D,bge}_{i,t} && \forall t \in \mathcal{T}\label{DSOp2pk}
\end{align}
\end{subequations}

\subsection{Bilevel Formulation}
Building upon the theoretical framework established in Section \ref{sec:Section_3}, the coordination between DSOs and the TSO is formalized as a bilevel Stackelberg problem, where the aggregated DSO layer acts as the upper-level leader and the TSO operates as the lower-level follower. This structure captures the hierarchical nature of decision-making while preserving the separation between local and system objectives described in Sections \ref{Subsec:TN_C} and \ref{Subsec:DN_C}. However, directly solving this bilevel model is computationally challenging due to the nested dependence between both levels. Thus, to obtain a tractable formulation, the lower-level problem of the TSO is substituted by its KKT optimality conditions. This reformulation relies on the convexity and continuity of the TSO’s problem and on the existence of a centralized coordination mechanism among DSOs, ensuring that the follower’s problem satisfies strong duality. In particular, since the TSO-level problem corresponds to a classical DC-OPF with quadratic generation costs and linear network constraints, 
the convexity and strong duality assumptions required for the KKT reformulation are satisfied~\cite{1216140}. Under these assumptions, the hierarchical structure presented in Section \ref{sec:Section_3} can be equivalently expressed as a single-level mathematical program, where the KKT conditions, comprising primal feasibility, dual feasibility, stationarity, and complementarity slackness, explicitly represent the TSO’s reaction to the upper-level decisions.

Therefore, the first component of the KKT system is the set of primal feasibility conditions, which reproduce the TSO's operational constraints, namely power balance, line flow physics, voltage and generation bounds, and transfer capacity at the TN-ADNs interfaces. These constraints are presented in~\eqref{FC} and are associated with Lagrange multipliers \( \lambda^T \) (for equality constraints) and dual variables \( \mu^T \) (for inequality constraints). These multipliers play a central role in forming the dual feasibility and stationarity conditions that complete the KKT system.

\begin{subequations}\label{FC}
\begin{flalign} 
    &\displaystyle\sum_{(i,j) \in \mathcal{L}_T} p^T_{i,j,t} - pg^T_{i,t} - pv^T_{i,t} + PL^T_{i,t} = 0 &&:\lambda^T_{1,i,t} \in \mathbb{R} &&& \forall i \in \Omega_T, \forall t \in \mathcal{T} \label{FCa}\\
    &\displaystyle\sum_{(i,j) \in \mathcal{L}_T} p^T_{i,j,t}  - pk^{T,bg}_{i,t} + pk^{T,sgc}_{i,t} + pk^{T,sge}_{i,t} = 0 && :\lambda^T_{2,i,t} \in \mathbb{R} &&& \forall i \in \Omega^T_{B}, \forall t \in \mathcal{T} \label{FCb}\\
    &X^T_{i,j} p^T_{i,j,t} - \theta^T_{i,t} + \theta^T_{j,t} = 0 &&:\lambda^T_{3,i,j,t} \in \mathbb{R} &&& \forall (i,j) \in \mathcal{L}_T, \, \forall t \in \mathcal{T} \label{FCc}\\
    & \theta^T_{i,t} = 0 &&:\lambda^T_{4,i,t} \in \mathbb{R} &&& \forall i \in \Omega^{T}_{ref}, \, \forall t \in \mathcal{T} \label{FCd}\\
    &p^T_{i,j,t} - S^{T,\max}_{i,j} \leq 0 &&: \mu^T_{1,i,j,t} \geq 0 &&& \forall (i,j) \in \mathcal{L}_T, \, \forall t \in \mathcal{T} \label{FCe}\\
    &- p^T_{i,j,t} - S^{T,\max}_{i,j} \leq 0 &&: \mu^T_{2,i,j,t} \geq 0 &&& \forall (i,j) \in \mathcal{L}_T, \, \forall t \in \mathcal{T} \label{FCf}\\
    &pg^T_{i,t} - PG^{T,\max}_{i,t} \leq 0 &&: \mu^T_{3,i,t} \geq 0 &&& \forall i \in \Omega_T, \forall t \in \mathcal{T} \label{FCg}\\
    &pv^T_{i,t} - PV^{T,\max}_{t} \gamma^{T,pv}_{i} \leq 0 &&: \mu^T_{4,i,t} \geq 0 &&& \forall i \in \Omega_T, \, \forall t \in \mathcal{T} \label{FCh}\\
    &\displaystyle\sum_{i \in \Omega^T_{B}}pk^{T,sgc}_{i,t} - \displaystyle\sum_{i \in \Omega_T} pv^{T}_{i,t} \leq 0  &&: \mu^T_{5,t} \geq 0 &&& \forall t \in \mathcal{T} \label{FCi}\\
    &\displaystyle\sum_{i \in \Omega^T_{B}}pk^{T,sge}_{i,t} - \displaystyle\sum_{i \in \Omega_T} pg^{T}_{i,t} \leq 0 &&: \mu^T_{6,t} \geq 0 &&& \forall t \in \mathcal{T} \label{FCj}\\
    & pk^{T,bg}_{i,t} - K^{T,bg}_{i} \leq 0 &&: \mu^T_{7,i,t} \geq 0 &&& \forall i \in \Omega^T_{B}, \forall t \in \mathcal{T} \label{FCk}
\end{flalign}
\end{subequations}

The second component of the KKT system is the stationarity conditions~\eqref{SC}, which enforce the vanishing gradient of the Lagrangian with respect to each decision variable. These conditions establish coupling between primal and dual variables and represent first-order optimality.
\begin{subequations}\label{SC}
\begin{flalign}
    &2Ca^T_{i} pg^T_{i,t} + Cb^T_{i} - \lambda^T_{1,i,t} + \mu^T_{3,i,t} - \mu^T_{6,t} = 0 && : pg^T_{i,t} &&& \forall i \in \Omega_T, \forall t \in \mathcal{T} \label{SCa}\\
    &\pi^T_{i} - \lambda^T_{1,i,t} + \mu^T_{4,i,t} - \mu^T_{5,t} = 0 && : pv^T_{i,t} &&& \forall i \in \Omega_T, \forall t \in \mathcal{T} \label{SCb}\\
    &\lambda^T_{2,i,t} + \mu^T_{5,t} = 0 && : pk^{T,sgc}_{i,t} &&& \forall i \in \Omega^T_{B}, \forall t \in \mathcal{T} \label{SCc}\\
    &\lambda^T_{2,i,t} + \mu^T_{6,t} = 0 && : pk^{T,sge}_{i,t} &&& \forall i \in \Omega^T_{B}, \forall t \in \mathcal{T} \label{SCd}\\
    &\mu^T_{7,i,t} - \lambda^T_{2,i,t} = 0 && : pk^{T,bg}_{i,t} &&& \forall i \in \Omega^T_{B}, \forall t \in \mathcal{T} \label{SCe}\\
    &\lambda^T_{1,i,t} + \lambda^T_{2,i,t} + \lambda^T_{3,i,j,t} X^T_{i,j} + \mu^T_{1,i,j,t} - \mu^T_{2,i,j,t}= 0 && : p^T_{i,j,t} &&& \forall (i,j) \in \mathcal{L}_T, \, \forall t \in \mathcal{T} \label{SCf}\\
    &\displaystyle\sum_{(i,j) \in \mathcal{L}_T} (\lambda^T_{3,j,i,t} - \lambda^T_{3,i,j,t}) + \lambda^T_{4,i,t} = 0 && :\theta^T_{i,t} &&& \forall i \in \Omega_T, \, \forall t \in \mathcal{T} \label{SCg}
\end{flalign}
\end{subequations}

Complementarity conditions, presented in~\eqref{BMC}, ensure that either the primal constraint is active or the corresponding dual variable is zero. Since these conditions are nonlinear and disjunctive, they are linearized using the Big-M method. This approach introduces binary variables \( \alpha^T \) and a sufficiently large constant \(  \mathcal{M}^T \) to capture the logic of disjunction. Although this increases the problem size and leads to a mixed-integer formulation, it retains the bilevel semantics in a form suitable for MILP solvers. The value of \( \mathcal{M}^T \) must be selected large enough to deactivate inactive constraints, yet not excessively large to avoid numerical instability or loose relaxation.

\begin{subequations}\label{BMC}
\begin{flalign}
    &p^T_{i,j,t} - S^{T,\max}_{i,j} \leq \mathcal{M}^T (1 - \alpha^T_{1,i,j,t})  && : \alpha^T_{1,i,j,t} \in \{0,1\} &&& \forall (i,j) \in \mathcal{L}_T, \, \forall t \in \mathcal{T} \label{BMCa}\\
    &\mu^T_{1,i,j,t} \leq \mathcal{M}^T (\alpha^T_{1,i,j,t})  && : \alpha^T_{1,i,j,t} \in \{0,1\} &&& \forall (i,j) \in \mathcal{L}_T, \, \forall t \in \mathcal{T} \label{BMCb}\\
    &- p^T_{i,j,t} - S^{T,\max}_{i,j} \leq \mathcal{M}^T (1 - \alpha^T_{2,i,j,t})  && : \alpha^T_{2,i,j,t} \in \{0,1\} &&& \forall (i,j) \in \mathcal{L}_T, \, \forall t \in \mathcal{T} \label{BMCc}\\
    &\mu^T_{2,i,j,t} \leq \mathcal{M}^T (\alpha^T_{2,i,j,t})  && : \alpha^T_{2,i,j,t} \in \{0,1\} &&& \forall (i,j) \in \mathcal{L}_T, \, \forall t \in \mathcal{T} \label{BMCd}\\
    &pg^T_{i,t} - PG^{T,\max}_{i,t} \leq \mathcal{M}^T (1 - \alpha^T_{3,i,t}) &&: \alpha^T_{3,i,t} \in \{0,1\} &&& \forall i \in \Omega_T, \forall t \in \mathcal{T} \label{BMCe}\\
    &\mu^T_{3,i,t} \leq \mathcal{M}^T (\alpha^T_{3,i,t}) && : \alpha^T_{3,i,t} \in \{0,1\} &&& \forall i \in \Omega_T, \forall t \in \mathcal{T} \label{BMCf}\\
    &pv^T_{i,t} - PV^{T,\max}_{t} \gamma^{T,pv}_{i} \leq \mathcal{M}^T (1 - \alpha^T_{4,i,t}) &&: \alpha^T_{4,i,t} \in \{0,1\} &&& \forall i \in \Omega_T, \forall t \in \mathcal{T} \label{BMCg}\\
    &\mu^T_{4,i,t} \leq \mathcal{M}^T (\alpha^T_{4,i,t}) && : \alpha^T_{4,i,t} \in \{0,1\} &&& \forall i \in \Omega_T, \forall t \in \mathcal{T} \label{BMCh}\\
    &\displaystyle\sum_{i \in \Omega^T_{B}}pk^{T,sgc}_{i,t} - \displaystyle\sum_{i \in \Omega_T} pv^{T}_{i,t} \leq \mathcal{M}^T (1 - \alpha^T_{5,t}) && : \alpha^T_{5,t} \in \{0,1\} &&& \forall t \in \mathcal{T} \label{BMCi}\\
    &\mu^T_{5,t} \leq \mathcal{M}^T (\alpha^T_{5,t}) && : \alpha^T_{5,t} \in \{0,1\} &&& \forall t \in \mathcal{T} \label{BMCj}\\
    &\displaystyle\sum_{i \in \Omega^T_{B}}pk^{T,sge}_{i,t} - \displaystyle\sum_{i \in \Omega_T} pg^{T}_{i,t} \leq \mathcal{M}^T (1 - \alpha^T_{6,t}) && : \alpha^T_{6,t} \in \{0,1\} &&& \forall t \in \mathcal{T} \label{BMCk}\\
    &\mu^T_{6,t} \leq \mathcal{M}^T (\alpha^T_{6,t}) && : \alpha^T_{6,t} \in \{0,1\} &&& \forall t \in \mathcal{T} \label{BMCl}\\
    & pk^{T,bg}_{i,t} - K^{T,bg}_{i} \leq \mathcal{M}^T (1 - \alpha^T_{7,i,t}) && : \alpha^T_{7,i,t} \in \{0,1\} &&& \forall i \in \Omega^T_{B}, \forall t \in \mathcal{T} \label{BMCm}\\
    &\mu^T_{7,i,t} \leq \mathcal{M}^T (\alpha^T_{7,i,t}) && : \alpha^T_{7,i,t} \in \{0,1\} &&& \forall i \in \Omega^T_{B}, \forall t \in \mathcal{T} \label{BMCn}
\end{flalign}
\end{subequations}

To ensure consistency between the TSO and DSOs at the interface level, we enforce boundary coupling constraints~\eqref{Boundary}, which equate the aggregated power injected or withdrawn at each boundary bus from both subsystems. This guarantees that energy transactions between the TN and each ADN are physically balanced and synchronized across the optimization layers.
\begin{subequations}\label{Boundary}
\begin{flalign}
    &pk^{D,sg}_{i,t} = pk^{T,bg}_{i,t} && \forall i \in \Omega^D_{B_k}, \forall t \in \mathcal{T} \label{Boundarya}\\
    &pk^{D,bgc}_{i,t} = pk^{T,sgc}_{i,t} && \forall i \in \Omega^D_{B_k}, \forall t \in \mathcal{T} \label{Boundaryb}\\
    &pk^{D,bge}_{i,t} = pk^{T,sge}_{i,t} && \forall i \in \Omega^D_{B_k}, \forall t \in \mathcal{T} \label{Boundaryc}
\end{flalign}
\end{subequations}

With the full KKT reformulation of the TSO's problem incorporated, the final model becomes a single-level MISOCP formulation that embeds the lower-level optimal response into the upper-level DSO coordination problem. This structure allows joint optimization over all DSOs’ decisions while respecting the physical and economic response of the TN, and is solvable using standard MILP techniques.

\subsection{MISOCP Reformulation}
The full coordination problem involving multiple DSOs and a single TSO is reformulated as a single-level MISOCP model. This formulation integrates the DSO-specific second-order cone network constraints and operational objectives with the KKT-based reformulation of the TSO’s lower-level convex problem. The resulting unified model simultaneously captures hierarchical interactions, resource constraints, and market mechanisms across the transmission and distribution layers. Therefore, the complete MISOCP formulation is expressed as:
\begin{align}
    \min \quad & 
    \displaystyle\sum_{t \in \mathcal{T}} \displaystyle\sum_{k\in \Omega^T_{B}} \displaystyle\sum_{i\in \Omega_{D_k}} (\lambda^{D,bgc}_{t} pk^{D,bgc}_{i,t} + \lambda^{D,bge}_{t}  pk^{D,bge}_{i,t} - \lambda^{D,sg}_{t} pk^{D,sg}_{i,t})
    \label{eq:final_milp} \\
    \text{s.t.} \quad 
    & \text{DSO Network Constraints (SOC-AC-OPF): } \eqref{DSOa}-\eqref{DSOi} \notag \\
    & \text{DSO Photovoltaic Generation Constraint: } \eqref{DSOpv} \notag \\
    & \text{DSO Peer-to-Peer Trading Constraints: } \eqref{DSOp2pa}-\eqref{DSOp2pk} \notag \\
    & \text{DSO Battery Operation Constraints: } \eqref{DSObta}-\eqref{DSObte} \notag \\
    & \text{TSO Primal Feasibility Conditions: } \eqref{FCa}-\eqref{FCk} \notag \\
    & \text{TSO Stationarity Conditions: } \eqref{SCa}-\eqref{SCg} \notag \\
    & \text{TSO Big-M Slack Conditions: } \eqref{BMCa}-\eqref{BMCn} \notag \\
    & \text{TSO-DSO Boundary Coupling Constraints: } \eqref{Boundarya}-\eqref{Boundaryc} \notag
\end{align}

\section{Case study and computational results}\label{sec:Section_5}

This section presents a case study that evaluates the proposed bilevel coordination framework under multiple operational configurations. The analysis begins by describing the hybrid TN–ADN system and then compares the coordination outcomes when the decision-making sequence is led either by the TSO or by the DSOs. Subsequently, the study examines the implicit competition among ADNs for low-cost generation available in the TN and assesses the influence of P2P energy trading mechanisms on network exchanges and active power flow reduction. Finally, the computational performance and scalability of the model are analyzed to show its tractability when the number of ADNs increases.

\subsection{Case study}
\begin{figure}[h]
    \centering
    \includegraphics[width=0.7\linewidth]{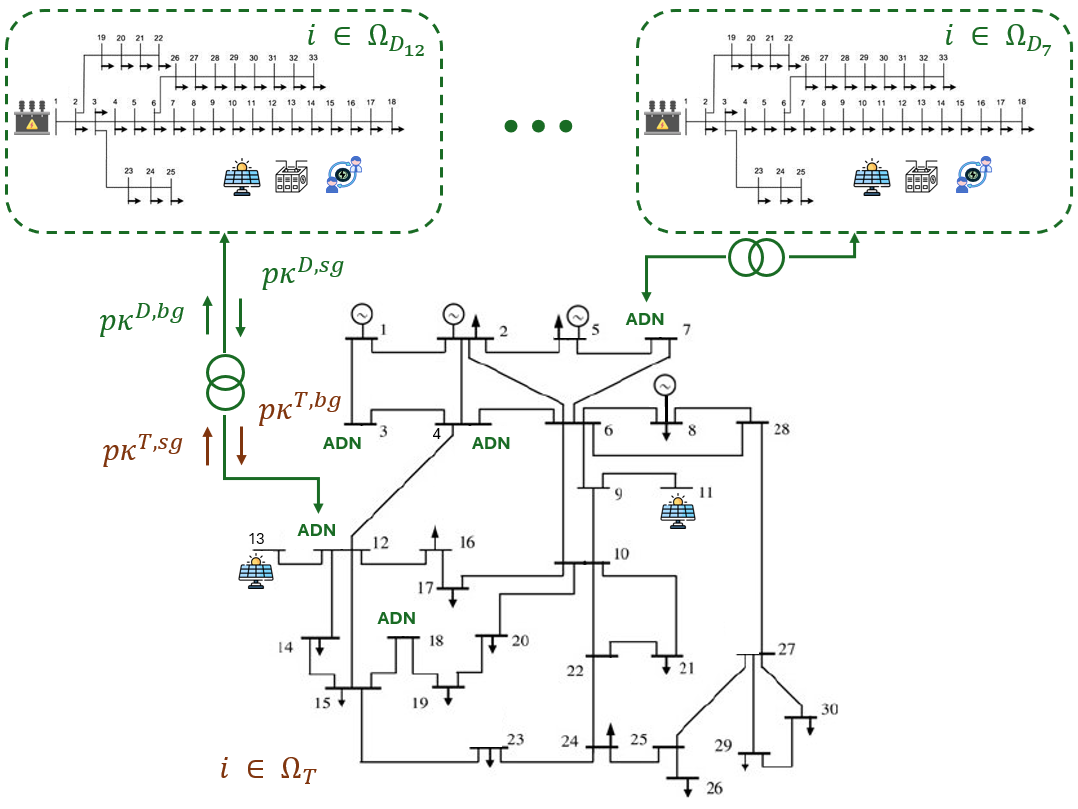}
    \caption{Schematic representation of the test system.}
    \label{fig:test_system_architecture_2}
\end{figure}

The proposed coordination model is evaluated using the IEEE 30-bus test system \cite{shahidehpour2003ieee} to represent the TN and the IEEE 33-bus test system \cite{9258930} to represent multiple ADNs. In this configuration, the TN includes controllable thermal generation units and large-scale PV plants, while each DN operates as an independent distribution area managed by a DSO considered as ADN with local PV and BESS resources. The models were executed over a 24-hour horizon representative of a day-ahead operation, with an hourly time step. In the TN, buses 3, 4, 7, 12, and 18 were selected as connection points for the different ADNs. Additionally, conventional generation units located at buses 11 and 13 of the transmission system were replaced by large-scale PV plants with rated capacities of 30 MW and 40 MW, respectively. The remaining demand buses of the TN were assumed to correspond to traditional distribution areas without DERs, thus behaving purely as passive loads. Consequently, these areas do not participate in the hierarchical decision-making process and maintain fixed demand profiles throughout the simulation. Figure~\ref{fig:test_system_architecture_2} depicts the hierarchical coupling between networks, highlighting the bidirectional active power exchanges at the interface buses through the boundary variables. The diagram distinguishes the elements and variables of each layer using a consistent color scheme: green represents the ADNs, including their internal components and exchanged variables, whereas orange denotes the TN and its associated variables. 

\begin{figure}[h]
    \centering
    \includegraphics[width=0.5\linewidth]{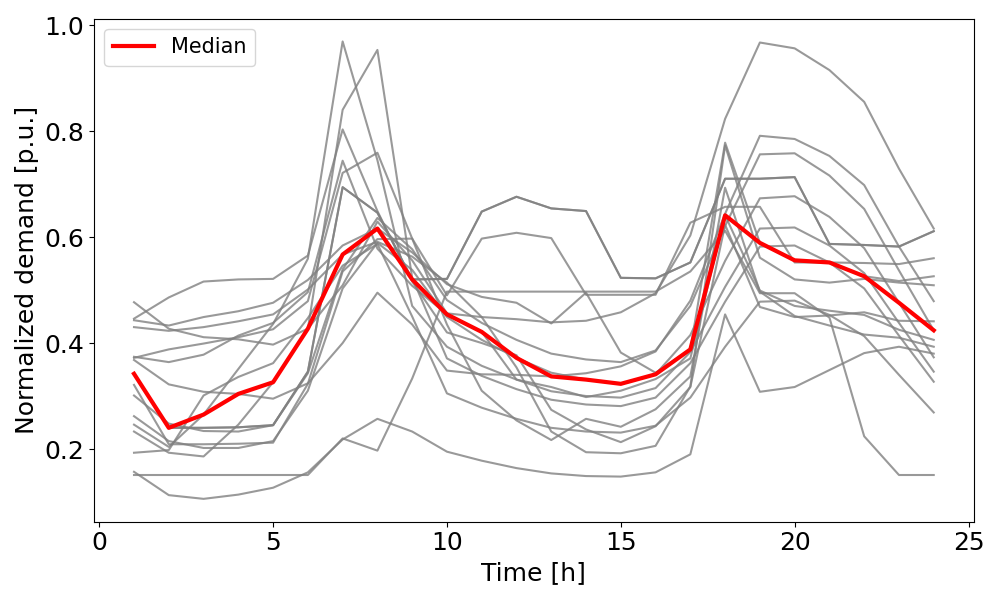}
    \caption{Normalized demand profiles}
    \label{fig:Demand_profiles}
\end{figure}

Seventeen normalized demand profiles, shown in Figure~\ref{fig:Demand_profiles}, were used to represent the electricity consumption patterns of both the TN and the ADNs. Each profile was scaled by the corresponding peak load of the associated network type, ensuring consistency between the TN and ADN demand levels. Since each ADN consists of 32 buses, the normalized profiles were cyclically assigned and repeated across all ADN nodes until completing the total number of demand points. 

Tables~\ref{TN_Data_G} and~\ref{DN_Data_G} summarize the main system characteristics, including peak load, generation capacity, and DER penetration levels across the networks. The aggregated peak demand of the five ADNs represents approximately 6.96\% of the total transmission-level demand, reflecting their contribution to the overall system load. At the TN level, the total PV installed capacity corresponds to 17\% of the system peak demand, illustrating a moderate renewable penetration scenario. Meanwhile, Table~\ref{DN_Data_G} shows that PV penetration among the ADNs varies from 35\% to 124\% when measured with respect to each ADN’s peak demand, capturing the heterogeneity of distributed generation across the ADNs.

\begin{table}[h]\centering
\begin{tabular}{lc}
\hline
\textbf{Item}             & \textbf{Q} \\ \hline
Nº Buses                     & 30         \\
Nº Branches                  & 41         \\
Peak Load system {[}MW{]} & 229.42     \\
Thermal gen {[}MW{]}      & 365        \\
PV gen {[}MW{]}           & 70         \\
Peak TD {[}MW{]}          & 213.46     \\
Peak ADN {[}MW{]}         & 15.96      \\
Nº Buses TD               & 20         \\
Nº buses con ADN               & 5          \\ \hline
\end{tabular}\caption{Transmission network general data. TD: traditional demand}\label{TN_Data_G}
\end{table}

\begin{table}[h]\centering
\begin{tabular}{lccccc}
\hline
                   & \textbf{DN 3} & \textbf{DN 4} & \textbf{DN 7} & \textbf{DN 12} & \textbf{DN 18} \\ \hline
\textbf{Peak Load [MW]} & 2.28          & 3.72          & 3.53          & 3.42           & 3.01           \\
\textbf{PV [MW]}        & 1.02          & 2             & 4.39          & 2.8            & 1.04           \\
\textbf{BESS [MWh]}      & 0.816         & 2             & 4.39          & 2.74           & 1.04           \\
\textbf{PV [\%]}     & 45\%          & 54\%          & 124\%         & 82\%           & 35\%           \\
\textbf{BESS [\%]}   & 36\%          & 54\%          & 124\%         & 80\%           & 35\%           \\ \hline
\end{tabular}\caption{Distribution networks general data}\label{DN_Data_G}
\end{table}

The model was implemented in Python using the Pyomo optimization framework, and all optimization problems were solved with Gurobi 12.0. Numerical simulations were executed on a workstation equipped with an Intel Core i9-14900HX processor (2.20 GHz) and 64 GB of RAM.

\subsection{DSO-led versus TSO-led}
\begin{figure}[htbp]
    \centering
    \includegraphics[width=0.6\linewidth]{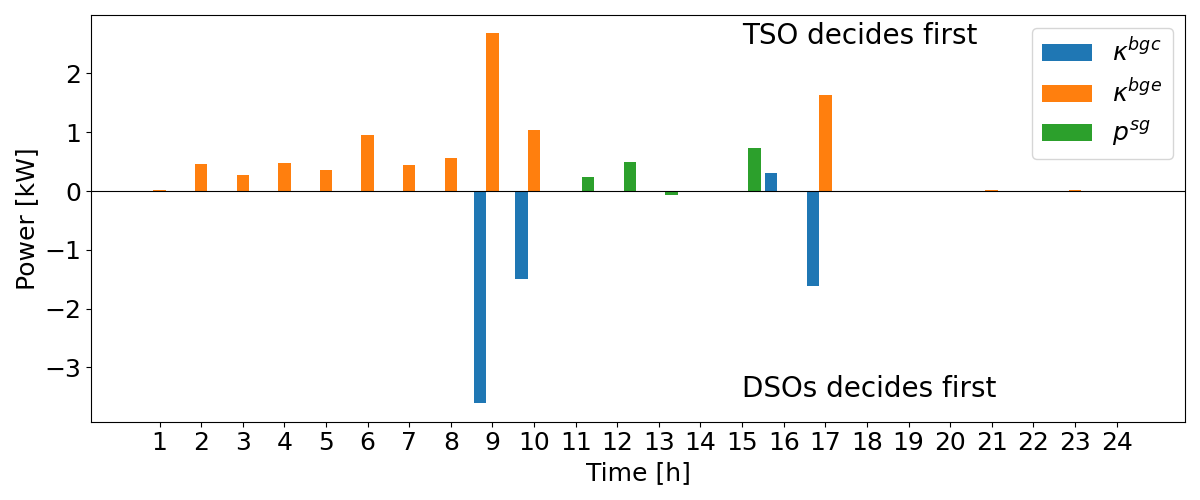}
    \caption{Comparison between the TSO-first and DSO-first decision sequences.}
    \label{fig:Kappas}
\end{figure}

This section compares the proposed hierarchical model, in which the DSOs act as leaders and the TSO as the follower, with the traditional top-down structure where the TSO decides first and the DSOs optimize based on that decision. However, formulating the top-down bilevel model would require relaxing the binary variables associated with the BESS operation and P2P trading within each ADN, resulting in a relaxed problem that is not directly comparable with the proposed DSO-leader formulation. Therefore, an alternative experiment was carried out. When the TSO decides first, each TN bus connected to an ADN is represented as a single aggregated node that combines the total demand of all community users and the installed PV and BESS capacities, excluding internal P2P exchanges. The TSO’s optimization problem thus assumes rational use of these aggregated DERs to minimize its economic dispatch, producing the values of $\kappa^{bg}$ and $p^{sg}$ for each ADN. These values are then fixed as parameters in the DSO problems, which are solved individually. The resulting solutions are compared with the DSO-leader case in Figure \ref{fig:Kappas}.

\begin{figure}[htbp]
    \centering
    \includegraphics[width=1.0\linewidth]{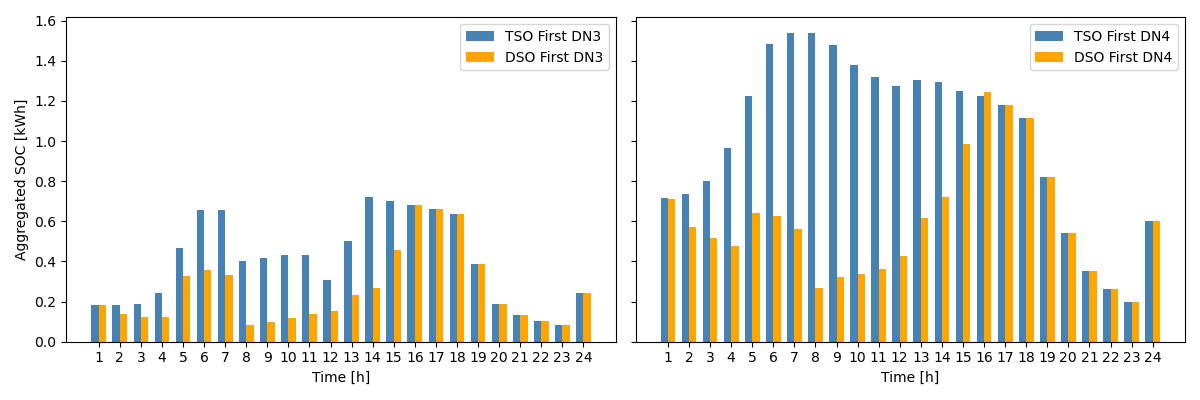}
    \caption{BESS ADN comparison between the TSO-first and DSO-first decision sequences}
    \label{fig:DERs}
\end{figure}

As shown in Figure \ref{fig:Kappas}, when the TSO decides first, the ADNs purchase more expensive energy $\kappa^{bge}$ from the TN and sell slightly more surplus 
$p^{sg}$, whereas when the DSOs decide first, they procure a higher share of cheaper energy $\kappa^{bgc}$. These differences in traded energy volumes result in an average increase of approximately 4.2\% in DSO operational costs when the TSO decides first. A similar experiment was conducted assuming equal purchase prices for dispatchable and PV generation, consistent with uniform day-ahead market clearing. Under this price assumption, the difference between $\kappa^{bg}$ components decreases, and the total cost gap between both hierarchical structures reduces to about 1\%. This occurs because, when the TSO decides first, it allocates the low-cost generation among all buses, including traditional loads in DNs and ADNs, whereas when the DSOs decide first, the cheaper resource is initially distributed among ADNs, and any surplus is then transferred to the remaining traditional DNs. Therefore, these results suggest that changing the decision order does not significantly affect the total economic dispatch cost. Nevertheless, Figure \ref{fig:DERs} highlights that the decision order does influence the operational pattern of BESS across the two hierarchical models.

Figure \ref{fig:DERs} shows the aggregated SoC of the BESS for two of the five ADN, connected to TN bus 3 (left) and bus 4 (right). The blue bars represent the SoC evolution when the TSO decides first, while the orange bars correspond to the DSO-first configuration. In both cases, the DSO-first sequence yields a more efficient use of storage: batteries charge earlier during low-cost hours and discharge closer to peak demand periods, improving local self-consumption and reducing dependence on expensive imports from the TN. Conversely, when the TSO decides first, the SoC trajectory is flatter and remains at higher levels throughout the day. This occurs because, in the TSO-first model, the operator assumes an aggregated demand profile together with aggregated DERs at each TN bus. As a result, the $\kappa^{bg}$ values obtained by the TSO do not reflect individual user decisions within the ADN. When these aggregated $\kappa^{bg}$ values are fixed in the DSO problems, they overestimate the required energy purchases between hours 2 and 16, forcing the BESS to store unnecessary surplus energy and maintain a higher SoC. This leads to a less adaptive and less efficient storage operation compared with the DSO-first configuration. This pattern is consistent across the remaining ADNs (see Table \ref{SOC_PL}), where the DSO-first configuration systematically achieves higher efficiency levels. For instance, in the ADN connected to bus 18, the aggregated SoC is 19.3\%  higher than the TSO-first sequence, reinforcing the idea that decentralized decision-making allows each DSO to better exploit local storage flexibility and adapt to its own demand and generation conditions.

\begin{table}[h]\centering
\begin{tabular}{cccccc}
\hline
                                 & \textbf{DN3} & \textbf{DN4} & \textbf{DN7} & \textbf{DN12} & \textbf{DN18} \\ \hline
TSO first               & 39.8\%       & 61.6\%       & 90.7\%       & 69.6\%        & 44.9\%        \\
DSO first               & 25.9\%       & 36.2\%       & 87.0\%       & 62.6\%        & 25.6\%        \\
$\Delta$ & 13.9\%       & 25.4\%       & 3.7\%        & 7.0\%         & 19.3\%        \\ \hline
\end{tabular}\caption{SoC regarding power load (PL)}\label{SOC_PL}
\end{table}

\subsection{Implicit competition}

This section analyzes the implicit competition that arises when multiple ADNs simultaneously demand the low-cost PV generation available at the TN level. To examine this effect, two experiments were performed. First, the hierarchical model was executed independently for each ADN (“Single” configuration), and then the complete model, including all five ADNs (“Full” configuration), was solved. Table \ref{tab:Comparison} summarizes the differences in energy purchased from the TN under both cases. The results show that the largest deviations occur in two time blocks that coincide with the beginning and end of PV generation hours. These periods correspond to the transition between low-cost PV and higher-cost thermal generation. When all ADNs are connected, they compete for access to the limited PV resource, which must be shared among them. As a result, each ADN obtains a smaller share of the cheaper energy, slightly reducing the economic benefit compared to the single-ADN case. In contrast, during mid-day hours (between 9:00 and 15:00), when PV production exceeds total demand, this competition effect disappears because the available cheap resource is sufficient to satisfy all DNs simultaneously.

\setlength{\tabcolsep}{2pt} 
\renewcommand{\arraystretch}{1.2} 
\begin{table}[]
\begin{tabular}{lccccccccccccccccccccccccc}
\hline
\textbf{} & \textbf{$DN^{3}$} & \textbf{1} & \textbf{2} & \textbf{3} & \textbf{4} & \textbf{5} & \textbf{6} & \textbf{7} & \textbf{8} & \textbf{9} & \textbf{10} & \textbf{11} & \textbf{12} & \textbf{13} & \textbf{14} & \textbf{15} & \textbf{16} & \textbf{17} & \textbf{18} & \textbf{19} & \textbf{20} & \textbf{21} & \textbf{22} & \textbf{23} & \textbf{24} \\ \hline
Single     & $\kappa^{sgc}$         & 0          & 0          & 0          & 0          & 0          & 0          & 0.7        & 1.5        & 1          & 0.6         & 0.3         & 0.3         & 0.1         & 0.1         & 0.2         & 0.5         & 0.9         & 1.4         & 0           & 0           & 0           & 0           & 0           & 0           \\
          & $\kappa^{sge}$         & 0.6        & 0.7        & 0.7        & 0.7        & 0.9        & 0.9        & 0.4        & 0          & 0          & 0           & 0           & 0           & 0           & 0           & 0           & 0           & 0           & 0           & 1.2         & 1.2         & 1.3         & 1.2         & 1.1         & 1           \\
          & $p^{bg}$          & 0          & 0          & 0          & 0          & 0          & 0          & 0          & 0          & 0          & 0           & 0           & 0           & 0           & 0           & 0           & 0           & 0           & 0           & 0           & 0           & 0           & 0           & 0           & 0           \\ \hline
Full  & $\kappa^{sgc}$         & 0          & 0          & 0          & 0          & 0          & 0          & 0          & 0          & 1          & 0.6         & 0.3         & 0.3         & 0.2         & 0.1         & 0.4         & 0.6         & 0.9         & 0           & 0           & 0           & 0           & 0           & 0           & 0           \\
          & $\kappa^{sge}$         & 0.6        & 0.7        & 0.7        & 0.8        & 1          & 0.9        & 1.2        & 1          & 0          & 0           & 0           & 0           & 0           & 0           & 0           & 0           & 0           & 1.4         & 1.2         & 1.2         & 1.3         & 1.2         & 1.1         & 1           \\
          & $p^{bg}$          & 0          & 0          & 0          & 0          & 0          & 0          & 0          & 0          & 0          & 0           & 0           & 0           & 0           & 0           & 0           & 0           & 0           & 0           & 0           & 0           & 0           & 0           & 0           & 0           \\ \hline
$\Delta$     & $\kappa^{sgc}$         & 0          & 0          & 0          & 0          & 0          & 0          & \cellcolor{red!20}-0.7       & \cellcolor{red!40}-1.5       & 0          & 0           & 0           & 0           & \cellcolor{myblue!05}0.1         & 0           & \cellcolor{myblue!05}0.1         & 0           & 0           & \cellcolor{red!40}-1.4        & 0           & 0           & 0           & 0           & 0           & 0           \\
          & $\kappa^{sge}$         & 0          & 0          & 0          & 0          & \cellcolor{myblue!05}0.1        & 0          & \cellcolor{myblue!20}0.8        & 1          & 0          & 0           & 0           & 0           & 0           & 0           & 0           & 0           & 0           & \cellcolor{myblue!40}1.4         & 0           & 0           & 0           & 0           & 0           & 0           \\ \hline
\end{tabular}\caption{Comparison of purchased energy between single- and full-ADN configurations} \label{tab:Comparison}
\end{table}

\begin{figure}[htbp]
    \centering
    \includegraphics[width=0.5\linewidth]{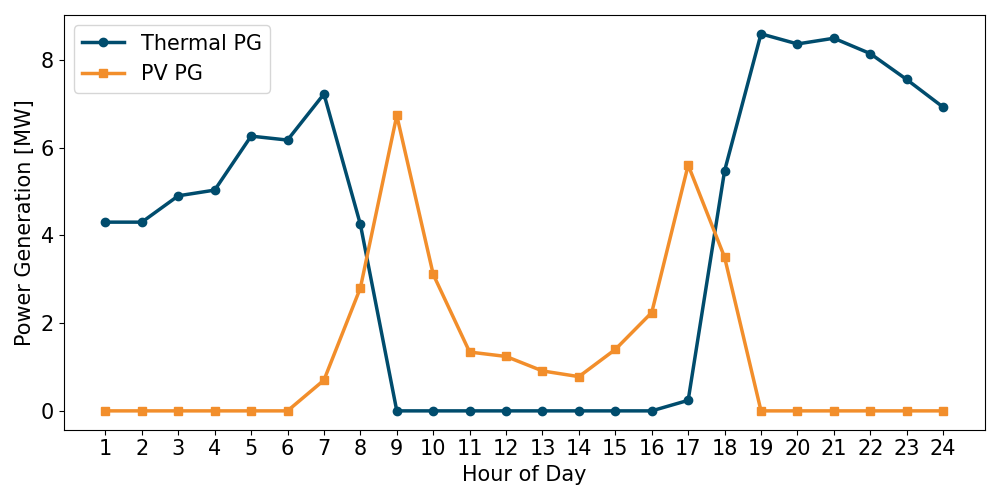}
    \caption{Hours in which PV and PG share a generation.}
    \label{fig:PG_sources}
\end{figure}

This interpretation is further supported by Figure \ref{fig:PG_sources}, which shows the aggregated energy purchased by the DSOs from the two generation sources available at the TN level. The figure clearly illustrates that the competition periods coincide with the hours when both PV and thermal generation are active. Once PV generation becomes available, thermal dispatch decreases sharply, and during mid-day hours, it drops to zero, as DSOs fully shift their purchases toward the cheaper PV source. The valley observed in the PV purchase curve at the TN level reflects the contribution of local PV generation within the ADNs, which reduces their need to buy energy from the upstream network.

\subsection{P2P energy trading and congestion effect}
This subsection examines several aspects related to the interaction between P2P energy trading and network operation within the proposed hierarchical framework. First, we analyze whether the decision sequence, whether the TSO or the DSOs act first, affects the total amount of energy traded within the ADNs. Second, we identify the specific time periods during inter-ADN exchanges, highlighting the temporal patterns of local and cross-network transactions. Finally, we assess how the presence of ADNs influences the active power flows across the transmission network by comparing the resulting dispatch with the baseline case in which the distribution areas operate under traditional consumption patterns.

Table \ref{P2P_PL} presents the percentage of energy traded within each ADN relative to its total uncontrollable demand. On average, the amount of energy exchanged through P2P transactions represents approximately 5.5\% of the total energy consumed within these ADNs. It can also be observed that when the TSO decides first, the volume of energy traded is higher than when the DSO acts first. This behavior is consistent with the discussion in the previous subsection: when the TSO leads the decision process, the energy surpluses committed at the transmission level must be stored in the BESS of each ADN. Consequently, the higher stored energy levels allow prosumers to offer more energy to their peers in subsequent hours. However, this also reveals a potential trade-off between promoting P2P exchanges (which in this case are indirectly driven by forced storage) and achieving an efficient use of distributed storage resources.

\begin{table}[h]\centering
\begin{tabular}{cccccc}
\hline
\textbf{} & \textbf{DN3} & \textbf{DN4} & \textbf{DN7} & \textbf{DN12} & \textbf{DN18} \\ \hline
TSO first & 7.9\%        & 8.2\%        & 6.6\%        & 4.1\%         & 5.2\%         \\
DSO first & 7.4\%        & 6.1\%        & 6.5\%        & 4.2\%         & 3.5\%         \\
$\Delta$       & 0.5\%        & 2.1\%        & 0.1\%        & -0.1\%        & 1.7\%         \\ \hline
\end{tabular}\caption{P2P energy trading regarding PL}\label{P2P_PL}
\end{table}

\begin{figure}[h]
    \centering
    \includegraphics[width=0.5\linewidth]{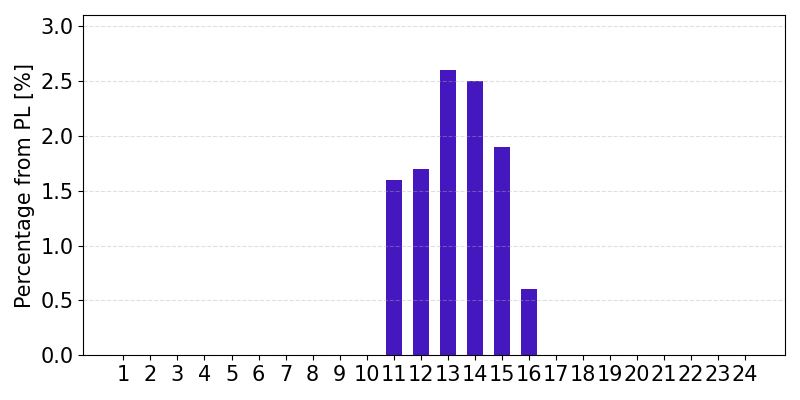}
    \caption{P2P energy trading at TN level}
    \label{fig:P2P_TN}
\end{figure}

Regarding the hourly distribution of inter-ADN trading, the bar chart in Figure \ref{fig:P2P_TN} shows that these exchanges primarily occur during solar generation hours, when the ADNs exhibit energy surpluses that can be traded with other DNs operating under traditional consumption patterns. The traded energy accounts for approximately 1\%–2.6\% of the total system demand (including both traditional loads and ADN demand), which is a non-negligible contribution considering that the aggregated peak load of the ADNs represents only 6.96\% of the total peak demand. This suggests that as the number of ADNs increases, the overall share of P2P energy trading in the system could become significantly higher.

To assess whether the integration of ADNs has a measurable effect on active power flows in the HV-TN, the hierarchical model was executed under two configurations: first, with all DNs operating under a traditional structure without PV and BESS resources; and second, with the inclusion of five ADNs equipped with distributed generation and storage. This comparison allows identifying how local DER operation modifies the power exchanges at the transmission level and whether it contributes to mitigating congestion in specific lines. Thus, Table \ref{COngestion_Reduction} summarizes a subset of transmission lines where the most significant flow reductions were observed after the transition from conventional DNs to ADNs. The results show that the line connecting buses 13–12 exhibits a decrease between 3\% and 6\% in active power flow during PV generation hours, while the line 9–11 shows a reduction of up to 4\%. These effects occur precisely when the ADNs achieve greater autonomy, since during hours of local PV generation they supply part of their demand internally and even inject surplus energy into the grid, in contrast to the case without DERs, where the same areas behaved purely as consumers.

\setlength{\tabcolsep}{4pt}
\begin{table}[h]
\begin{tabular}{cccccccccccccccccccccccccc}
\hline
\textbf{from} & \textbf{to} & \textbf{1} & \textbf{2} & \textbf{3} & \textbf{4} & \textbf{5} & \textbf{6} & \textbf{7} & \textbf{8} & \textbf{9} & \textbf{10} & \textbf{11} & \textbf{12} & \textbf{13} & \textbf{14} & \textbf{15} & \textbf{16} & \textbf{17} & \textbf{18} & \textbf{19} & \textbf{20} & \textbf{21} & \textbf{22} & \textbf{23} & \textbf{24} \\ \hline
\textbf{13}   & \textbf{12} & 0        & 0        & 0        & 0        & 0        & 0        & 0        & 0        & 0        & \cellcolor{myblue!30}3         & \cellcolor{myblue!60}6         & \cellcolor{myblue!40}4         & \cellcolor{myblue!50}5         & \cellcolor{myblue!40}4         & \cellcolor{myblue!40}4         & \cellcolor{myblue!30}3         & 0         & 0         & 0         & 0         & 0         & 0         & 0         & 0         \\
\textbf{11}   & \textbf{9}  & 0        & 0        & 0        & 0        & 0        & 0        & 0        & 0        & 0        & \cellcolor{myblue!30}3         & \cellcolor{myblue!40}4         & \cellcolor{myblue!30}3         & \cellcolor{myblue!40}4         & \cellcolor{myblue!30}3         & \cellcolor{myblue!30}3         & \cellcolor{myblue!30}3         & 0         & 0         & 0         & 0         & 0         & 0         & 0         & 0         \\
\textbf{15}   & \textbf{18} & \cellcolor{myblue!10}1        & 0        & 0        & 0        & 0        & 0        & 0        & \cellcolor{myblue!10}1        & 0        & \cellcolor{myblue!10}1         & \cellcolor{myblue!30}3         & \cellcolor{myblue!20}2         & \cellcolor{myblue!20}2         & \cellcolor{myblue!20}2         & \cellcolor{myblue!20}2         & \cellcolor{myblue!10}1         & 0         & 0         & \cellcolor{myblue!10}1         & \cellcolor{myblue!10}1         & 0         & 0         & 0         & \cellcolor{myblue!10}1         \\
\textbf{6}    & \textbf{9}  & 0        & 0        & 0        & 0        & 0        & 0        & 0        & 0        & 0        & \cellcolor{myblue!20}2         & \cellcolor{myblue!30}3         & \cellcolor{myblue!20}2         & \cellcolor{myblue!20}2         & \cellcolor{myblue!20}2         & \cellcolor{myblue!20}2         & \cellcolor{myblue!20}2         & 0         & 0         & 0         & 0         & 0         & 0         & 0         & 0         \\
\textbf{6}    & \textbf{7}  & 0        & 0        & 0        & 0        & 0        & 0        & 0        & 0        & 0        & \cellcolor{myblue!10}1         & \cellcolor{myblue!20}2         & \cellcolor{myblue!20}2         & \cellcolor{myblue!20}2         & \cellcolor{myblue!20}2         & \cellcolor{myblue!10}1         & \cellcolor{myblue!10}1         & 0         & 0         & 0         & 0         & 0         & 0         & 0         & 0         \\
\textbf{4}    & \textbf{12} & 0        & 0        & 0        & 0        & 0        & 0        & 0        & \cellcolor{myblue!10}1        & 0        & 0         & \cellcolor{myblue!20}2         & \cellcolor{myblue!10}1         & \cellcolor{myblue!10}1         & \cellcolor{myblue!10}1         & \cellcolor{myblue!20}2         & \cellcolor{myblue!10}1         & 0         & 0         & 0         & 0         & 0         & 0         & 0         & 0         \\
\textbf{1}    & \textbf{2}  & \cellcolor{myblue!10}1        & 0        & 0        & 0        & 0        & 0        & 0        & \cellcolor{myblue!10}1        & \cellcolor{myblue!10}1        & 0         & 0         & \cellcolor{myblue!10}1         & \cellcolor{myblue!10}1         & \cellcolor{myblue!10}1         & 0         & 0         & \cellcolor{myblue!10}1         & \cellcolor{myblue!10}1         & \cellcolor{myblue!10}1         & \cellcolor{myblue!10}1         & 0         & 0         & 0         & \cellcolor{myblue!10}1         \\
\textbf{9}    & \textbf{10} & 0        & 0        & 0        & 0        & 0        & 0        & 0        & 0        & 0        & \cellcolor{myblue!10}1         & \cellcolor{myblue!10}1         & \cellcolor{myblue!10}1         & \cellcolor{myblue!10}1         & \cellcolor{myblue!10}1         & \cellcolor{myblue!10}1         & \cellcolor{myblue!10}1         & 0         & 0         & 0         & 0         & 0         & 0         & 0         & 0         \\
\textbf{27}   & \textbf{25} & 0        & 0        & 0        & 0        & 0        & 0        & 0        & 0        & 0        & \cellcolor{myblue!10}1         & \cellcolor{myblue!20}2         & \cellcolor{myblue!10}1         & \cellcolor{myblue!10}1         & \cellcolor{myblue!10}1         & \cellcolor{myblue!10}1         & \cellcolor{myblue!10}1         & 0         & 0         & 0         & 0         & 0         & 0         & 0         & 0         \\
\textbf{10}   & \textbf{20} & 0        & 0        & 0        & 0        & 0        & 0        & 0        & 0        & 0        & \cellcolor{myblue!10}1         & \cellcolor{myblue!10}1         & \cellcolor{myblue!10}1         & \cellcolor{myblue!10}1         & \cellcolor{myblue!10}1         & \cellcolor{myblue!10}1         & 0         & 0         & 0         & 0         & 0         & 0         & 0         & 0         & 0         \\
\textbf{20}   & \textbf{19} & 0        & 0        & 0        & 0        & 0        & 0        & 0        & 0        & 0        & \cellcolor{myblue!10}1         & \cellcolor{myblue!10}1         & \cellcolor{myblue!10}1         & \cellcolor{myblue!10}1         & \cellcolor{myblue!10}1         & \cellcolor{myblue!10}1         & 0         & 0         & 0         & 0         & 0         & 0         & 0         & 0         & 0         \\
\textbf{22}   & \textbf{24} & 0        & 0        & 0        & 0        & 0        & 0        & 0        & 0        & 0        & \cellcolor{myblue!10}1         & \cellcolor{myblue!10}1         & \cellcolor{myblue!10}1         & \cellcolor{myblue!10}1         & \cellcolor{myblue!10}1         & \cellcolor{myblue!10}1         & \cellcolor{myblue!10}1         & 0         & 0         & 0         & 0         & 0         & 0         & 0         & 0   
\\ \hline
\end{tabular}\caption{Percentage reduction of active power flows in the TN caused by DER operation in ADNs} \label{COngestion_Reduction}
\end{table}

\begin{figure}[h]
    \centering
    \includegraphics[width=0.5\linewidth]{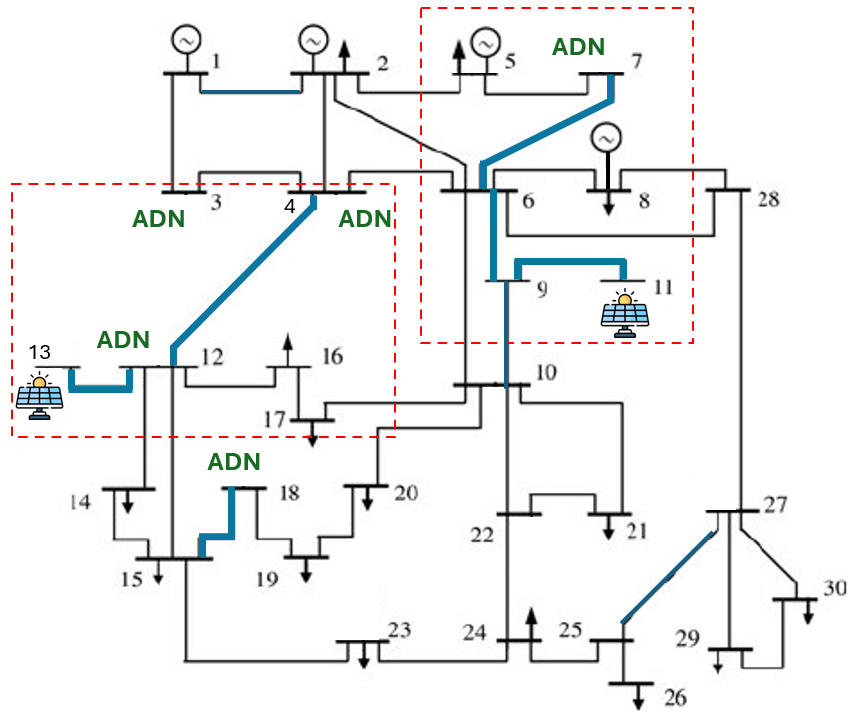}
    \caption{Congestion}
    \label{fig:Congestion}
\end{figure}

Figure \ref{fig:Congestion} highlights in blue the transmission lines affected by these reductions, grouped into two main clusters (red boxes) that correspond to the areas where ADNs and PV generation are located. The spatial concentration of these effects indicates that the presence of DERs not only alleviates local loading conditions but also modifies upstream power transfers in adjacent areas. Additionally, smaller reductions are visible in other lines outside the red clusters, suggesting that ADN operation induces a broader redistribution of active power across the TN, with potential benefits in congestion management and overall network efficiency.

\subsection{Model performance}
\begin{figure}[h]
    \centering
    \includegraphics[width=0.6\linewidth]{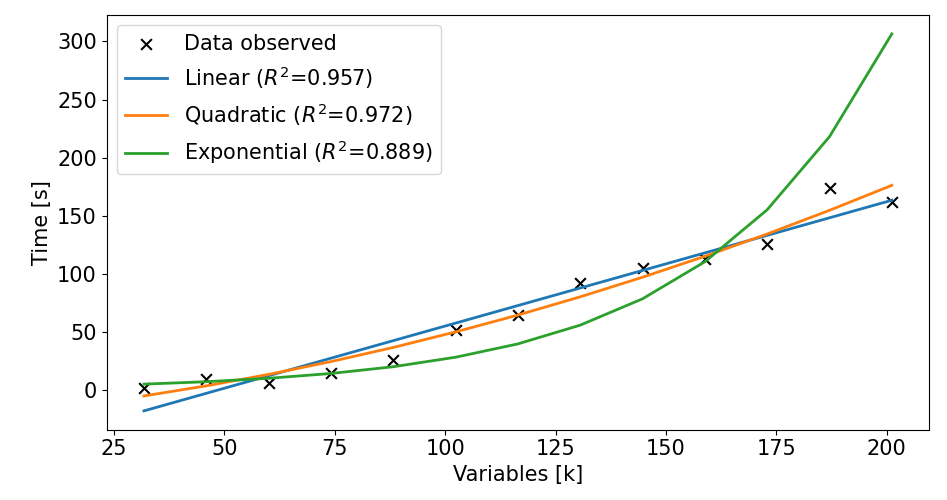}
    \caption{Model fitting to solution time}
    \label{fig:Plot_Times}
\end{figure}

This subsection extends the previous experiments by increasing the number of ADNs connected to the transmission grid to evaluate the scalability of the proposed coordination model. While the earlier analyses considered five ADNs, additional networks were progressively incorporated until reaching a configuration with thirteen ADNs. This expansion results in a hybrid system comprising 459 buses and more than 200,000 decision variables, including both continuous and binary ones. Thus, Table \ref{Perfomance_Table}  and Figure \ref{fig:Plot_Times} illustrate the computational performance of the model as the number of ADNs increases. The execution time is plotted as a function of the total number of variables, providing a reference for assessing the growth trend of computational complexity. The results show that the model’s runtime exhibits a near-quadratic, rather than exponential—growth pattern, remaining within tractable limits even as the system size expands. This behavior suggests that the proposed DSO-led bilevel formulation maintains favorable growing properties despite including discrete decisions related to DER coordination and P2P trading.

\begin{table}[h]
\centering
\begin{tabular}{ccccc}
\hline
\textbf{Nº ADNs   connected to TN} & \textbf{Total Buses} & \textbf{Variables} & \textbf{Constraints} & \textbf{Times {[}s{]}} \\ \hline
1                                  & 63                   & 31,800             & 39,312               & 2.1                    \\
2                                  & 96                   & 45,912             & 55,560               & 9.7                    \\
3                                  & 129                  & 60,024             & 71,856               & 6.2                    \\
4                                  & 162                  & 74,136             & 88,080               & 14.7                   \\
5                                  & 195                  & 88,248             & 104,232              & 25.7                   \\
6                                  & 228                  & 102,360            & 120,912              & 52.1                   \\
7                                  & 261                  & 116,472            & 137,592              & 64.9                   \\
8                                  & 294                  & 130,584            & 154,272              & 92.1                   \\
9                                  & 327                  & 144,696            & 170,952              & 104.9                  \\
10                                 & 360                  & 158,808            & 187,632              & 112.7                  \\
11                                 & 393                  & 172,920            & 204,312              & 125.3                  \\
12                                 & 426                  & 187,032            & 220,992              & 173.5                  \\
13                                 & 459                  & 201,144            & 237,192              & 162.2                  \\ \hline
\end{tabular}\caption{Resolution times involved when the number of ADN connected to the TN increases}\label{Perfomance_Table}
\end{table}

Although these results correspond to a specific case study and a given level of DER penetration, including PV, BESS, and P2P energy trading, they suggest that the complexity introduced by distribution-level binary variables can be effectively managed within the proposed modeling framework, without resorting to the variable relaxations typically required in conventional top-down (TSO-first) formulations. Consequently, the proposed approach achieves an accurate representation of ADN operational details while maintaining computational efficiency, thus representing a suitable scheme for large-scale multi-ADN systems.

\section{Conclusions}\label{sec:Section_6}
This work proposed a bi-level coordination model that captures the hierarchical interaction between the transmission and distribution layers in systems with multiple ADNs. The proposed framework represents a DSO-led coordination scheme, in which each DSO optimizes its local operation before the TSO performs the system coordination. The model integrates the physical constraints of both network layers and the P2P energy exchanges within and across ADNs through the TN, offering a computationally tractable formulation without relaxing the binary variables that define the operational states of DERs within the ADNs.

The case study showed that when the DSOs decide first, the batteries tend to be used more efficiently than under the traditional top-down coordination scheme, where the TSO decides first with limited information about ADN users. This behavior improves local self-consumption and reduces energy imports from the TN. Moreover, the introduction of ADNs could effectively mitigate congestion in the TN lines by reducing the net power exchanged at the interface buses, particularly during PV generation hours when some consumption points become net generators at the transmission level. The computational analysis also indicated that the model scales approximately linearly with the number of ADNs, maintaining solvability without relaxing binary variables, which suggests its potential applicability to larger interconnected systems.

Finally, future research could explore two main directions. First, although the proposed model adopts a centralized implementation of the DSO-led coordination, a decentralized resolution could be developed to parallelize the optimization problems of individual DSOs. Such an approach could improve computational efficiency, provided the hierarchical structure, where DSOs decide first and the TSO reacts, remains preserved. Second, the present deterministic formulation can be extended to incorporate uncertainty in distributed generation, demand, and emerging flexible assets such as electric vehicles. Including these stochastic elements would allow a more realistic assessment of operational and planning decisions under variable conditions in ADNs.

\section*{Acknowledgment}
This work has been supported by ANID FONDECYT Iniciación 11240745.

\begin{table*}[!t]
\footnotesize
\begin{framed}
\begin{longtable}{p{3cm}p{10cm}}
\multicolumn{2}{l}{\textbf{Sets}}\\[2pt]
$i \in \Omega_T$ & Set of buses in the TN.\\
$k \in \Omega^T_{B}$ & Set of boundary buses connecting the TN with each DN; $\Omega^T_{B} \subset \Omega_T$.\\
$i \in \Omega^{T}_{ref}$ & Set of reference buses in the TN; $\Omega^{T}_{ref} \subset \Omega_T$.\\
$i \in \Omega_{D_k}$ & Set of buses in the DN associated with boundary bus $k \in \Omega^T_{B}$.\\
$i \in \Omega^D_{A_k}$ & Set of agents in DN $k$; $\Omega^D_{A_k} \subseteq \Omega_{D_k}$.\\
$i \in \Omega^D_{B_k}$ & Set of boundary buses within DN $k$; $\Omega^D_{B_k} \subset \Omega_{D_k}$.\\
$(i,j) \in \mathcal{L}_T$ & Set of lines in the TN; $\mathcal{L}_T = \{(i,j)\mid i,j \in \Omega_T\}$.\\
$(i,j) \in \mathcal{L}_{D_k}$ & Set of lines in DN $k$; $\mathcal{L}_{D_k} = \{(i,j)\mid i,j \in \Omega_{D_k}\}$.\\
$t \in \mathcal{T}$ & Set of discrete time periods.\\
$\phi$ & Index representing the network layer, $\phi \in \{T,D\}$.\\[4pt]

\multicolumn{2}{l}{\textbf{Parameters}}\\[2pt]
$Ca^T_i, Cb^T_i, Cc^T_i$ & Cost coefficients of TN generation [\$/MW$^2$, \$/MW, \$].\\
$\pi^T_i$ & Marginal cost of PV generation at TN bus $i$ [\$/MWh].\\
$PG^{T,\max}_{i,t}$ & Max active power generation at TN bus $i$ [MW].\\
$S^{T,\max}_{i,j}$ & Max apparent power on TN line $(i,j)$ [MVA].\\
$X^{\phi}_{i,j}, R^D_{i,j}$ & Reactance and resistance of line $(i,j)$ [p.u.].\\
$QL^D_{i,t}$ & Reactive demand at DN bus $i$ [MVAr].\\
$QG^{D,\min}_{i}, QG^{D,\max}_{i}$ & Limits of reactive generation at DN bus $i$ [MVAr].\\
$\lambda^{D,sg}_t, \lambda^{D,bgc}_t, \lambda^{D,bge}_t$ & Energy prices in DN at time $t$ [\$/MWh].\\
$V^{D,\min}_{i}, V^{D,\max}_{i}$ & Voltage limits at DN bus $i$ [p.u.].\\
$PL^{\phi}_{i,t}$ & Active power demand at TN/DN bus $i$ [MW].\\
$PV^{\phi,\max}_{t}$ & Max PV generation available [MW].\\
$\gamma^{\phi,pv}_i$ & Installed PV capacity ratio at bus $i$.\\
$PB^{D,bt}$ & Max charge/discharge power of battery [MW].\\
$SOC^{D,bt}_{\min/\max}$ & Min/Max SOC of battery [p.u.].\\
$\gamma^{D,bt}_{i}$ & Battery capacity at agent $i$ [MWh].\\
$\nu^{D,bt}_{i}$ & Binary parameter: 1 if a battery is installed at $i$.\\
$\varphi^{D,ch}, \varphi^{D,ds}$ & Charging/discharging efficiencies.\\
$K^{T,bg}_{i}$ & Max power bought by TN at boundary bus $i$ [MW].\\
$K^{D,sg}_{i}, K^{D,bgc}_{i}, K^{D,bge}_{i}$ & Max power sold or bought by DN at boundary buses [MW].\\[4pt]

\multicolumn{2}{l}{\textbf{Variables}}\\[2pt]
$pg^T_{i,t}$ & Active generation at TN bus $i$ [MW].\\
$pv^{\phi}_{i,t}$ & PV generation at TN/DN bus $i$ [MW].\\
$\theta^T_{i,t}$ & Voltage angle at TN bus $i$ [rad].\\
$p^T_{i,j,t}$ & Active power flow on TN line $(i,j)$ [MW].\\
$p^{T,sgc}_{i,t}, p^{T,sge}_{i,t}, p^{T,bg}_{i,t}$ & Power sold/bought between TN and DNs [MW].\\
$q^D_{i,j,t}$ & Reactive flow on DN line $(i,j)$ [MVAr].\\
$\ell^D_{i,j,t}$ & Squared current on DN line $(i,j)$ [A$^2$].\\
$qg^D_{i,t}$ & Reactive generation at DN bus $i$ [MVAr].\\
$v^D_{i,t}$ & Voltage magnitude at DN bus $i$ [p.u.].\\
$soc^{D,bt}_{i,t}$ & State of charge of battery at agent $i$ [p.u.].\\
$w^{D,bt}_{i,t}$ & Binary: battery charge/discharge mode.\\
$ch^{D,bt}_{i,t}, ds^{D,bt}_{i,t}$ & Charging/discharging power [MW].\\
$\Delta p^D_{i,t}, \Delta p^{D,B}_{i,t}$ & Net power exchange by agent/boundary [MW].\\
$\Delta p^{D,+}_{i,t}, \Delta p^{D,-}_{i,t}$ & Power sold/purchased by agent [MW].\\
$p^{D,sg}_{i,t}, p^{D,bg}_{i,t}$ & Power sold/bought between agent and DN boundary [MW].\\
$p^{D,sm}_{i,t}, p^{D,bm}_{i,t}$ & Power sold/bought in local P2P market [MW].\\
$y^D_{i,t}$ & Binary: 1 if agent $i$ sells at $t$.\\
\end{longtable}
\end{framed}
\caption{Nomenclature.}\label{Nomenclature_Table}
\end{table*}

\bibliographystyle{elsarticle-num-names}
\bibliography{references} 

@article{PAPADIS2020118025,
title = {Challenges in the decarbonization of the energy sector},
journal = {Energy},
volume = {205},
pages = {118025},
year = {2020},
issn = {0360-5442},
doi = {https://doi.org/10.1016/j.energy.2020.118025},
url = {https://www.sciencedirect.com/science/article/pii/S0360544220311324},
author = {Elisa Papadis and George Tsatsaronis}
}

@misc{irena2024a,
  author       = {{International Renewable Energy Agency (IRENA)}},
  year         = {2024},
  title        = {Energy Transition Outlook},
  url          = {https://www.irena.org/Energy-Transition/Outlook},
  note         = {Retrieved November 18, 2024}
}

@techreport{lind2024tso,
  author      = {Leandro Lind},
  title       = {TSO-DSO Coordination: A Multidimensional Study on Coordination Schemes, Modelling and Regulation in the European Context},
  institution = {Universidad Pontificia Comillas},
  year        = {2024},
  type        = {Master's Thesis},
  url         = {https://repositorio.comillas.edu/xmlui/handle/11531/95589}
}

@INPROCEEDINGS{9941459,
  author={Papazoglou, Georgios K. and Bakirtzis, Emmanouil A. and Forouli, Aikaterini A. and Biskas, Pandelis N. and Bakirtzis, Anastasios G.},
  booktitle={2022 2nd International Conference on Energy Transition in the Mediterranean Area (SyNERGY MED)}, 
  title={A two-stage market-based TSO-DSO coordination framework}, 
  year={2022},
  volume={},
  number={},
  pages={1-6},
  doi={10.1109/SyNERGYMED55767.2022.9941459}}

@article{lind2019transmission,
  author  = {Lind, Leandro and Cossent, Rafael and Chaves-Ávila, José Pablo and Gómez San Román, Tomás},
  title   = {Transmission and distribution coordination in power systems with high shares of distributed energy resources providing balancing and congestion management services},
  journal = {Wiley Interdisciplinary Reviews: Energy and Environment},
  volume  = {8},
  number  = {6},
  pages   = {e357},
  year    = {2019},
  doi     = {10.1002/wene.357},
  publisher = {Wiley Online Library}
}

@Article{en15197312,
AUTHOR = {Alazemi, Talal and Darwish, Mohamed and Radi, Mohammed},
TITLE = {TSO/DSO Coordination for RES Integration: A Systematic Literature Review},
JOURNAL = {Energies},
VOLUME = {15},
YEAR = {2022},
NUMBER = {19},
ARTICLE-NUMBER = {7312},
URL = {https://www.mdpi.com/1996-1073/15/19/7312},
ISSN = {1996-1073},
DOI = {10.3390/en15197312}
}

@INPROCEEDINGS{10202840,
  author={Simoes, Micael and Madureira, Andre G. and Soares, Filipe and Lopes, Joao Pecas},
  booktitle={2023 IEEE Belgrade PowerTech}, 
  title={TSO-DSO Coordinated Operational Planning in the Presence of Shared Resources $i,m,o,t$}, 
  year={2023},
  volume={},
  number={},
  pages={01-08},
  doi={10.1109/PowerTech55446.2023.10202840}}

@INPROCEEDINGS{10609023,
  author={Bjerland, Siri and Del Granado, Pedro Crespo and GrØttum, Hanne and Nokandi, Ehsan},
  booktitle={2024 20th International Conference on the European Energy Market (EEM)}, 
  title={TSO-DSO Coordination Under Wind and Solar Power Uncertainty: A Two-Stage Stochastic Programming Approach}, 
  year={2024},
  volume={},
  number={},
  pages={1-8},
  doi={10.1109/EEM60825.2024.10609023}}

@article{YUAN2017600,
title = {Hierarchical coordination of TSO-DSO economic dispatch considering large-scale integration of distributed energy resources},
journal = {Applied Energy},
volume = {195},
pages = {600-615},
year = {2017},
issn = {0306-2619},
doi = {https://doi.org/10.1016/j.apenergy.2017.03.042},
url = {https://www.sciencedirect.com/science/article/pii/S030626191730274X},
author = {Zhao Yuan and Mohammad Reza Hesamzadeh}
}

@article{SOARES2020100333,
title = {Reactive power provision by the DSO to the TSO considering renewable energy sources uncertainty},
journal = {Sustainable Energy, Grids and Networks},
volume = {22},
pages = {100333},
year = {2020},
issn = {2352-4677},
doi = {https://doi.org/10.1016/j.segan.2020.100333},
url = {https://www.sciencedirect.com/science/article/pii/S2352467719304886},
author = {Tiago Soares and Leonel Carvalho and Hugo Morais and Ricardo J. Bessa and Tiago Abreu and Eric Lambert}
}

@article{BERALDOBANDEIRA2024110818,
title = {An ADP framework for flexibility and cost aggregation: Guarantees and open problems},
journal = {Electric Power Systems Research},
volume = {234},
pages = {110818},
year = {2024},
issn = {0378-7796},
doi = {https://doi.org/10.1016/j.epsr.2024.110818},
url = {https://www.sciencedirect.com/science/article/pii/S0378779624007041},
author = {Maísa {Beraldo Bandeira} and Timm Faulwasser and Alexander Engelmann}
}

@ARTICLE{9503337,
  author={Bragin, Mikhail A. and Dvorkin, Yury},
  journal={IEEE Transactions on Power Systems}, 
  title={TSO-DSO Operational Planning Coordination Through “$l_1-$Proximal” Surrogate Lagrangian Relaxation}, 
  year={2022},
  volume={37},
  number={2},
  pages={1274-1285},
  doi={10.1109/TPWRS.2021.3101220}}

@ARTICLE{9835136,
  author={Weng, Yu and Xie, Jiahang and Wang, Peng and Nguyen, Hung Dinh},
  journal={IEEE Transactions on Power Systems}, 
  title={Asymmetrically Reciprocal Effects and Congestion Management in TSO-DSO Coordination Through Feasibility Regularizer}, 
  year={2023},
  volume={38},
  number={2},
  pages={1948-1962},
  doi={10.1109/TPWRS.2022.3193052}}

@ARTICLE{9914683,
  author={Steriotis, Konstantinos and Makris, Prodromos and Tsaousoglou, Georgios and Efthymiopoulos, Nikolaos and Varvarigos, Emmanouel},
  journal={IEEE Transactions on Power Systems}, 
  title={Co-Optimization of Distributed Renewable Energy and Storage Investment Decisions in a TSO-DSO Coordination Framework}, 
  year={2023},
  volume={38},
  number={5},
  pages={4515-4529},
  doi={10.1109/TPWRS.2022.3212919}}

@INPROCEEDINGS{9543104,
  author={Olsen, Ole Kjærland and Sieraszewski, Damian and Ivanko, Dmytro and Oleinikova, Irina and Farahmand, Hossein},
  booktitle={2021 International Conference on Smart Energy Systems and Technologies (SEST)}, 
  title={Hybrid AC/DC Optimal Power Flow Modelling Approach for Coordination in Flexibility Market}, 
  year={2021},
  volume={},
  number={},
  pages={1-6},
  doi={10.1109/SEST50973.2021.9543104}}

@article{RODRIGUES2023101204,
title = {Reactive power management considering Transmission System Operator and Distribution System Operator coordination},
journal = {Sustainable Energy, Grids and Networks},
volume = {36},
pages = {101204},
year = {2023},
issn = {2352-4677},
doi = {https://doi.org/10.1016/j.segan.2023.101204},
url = {https://www.sciencedirect.com/science/article/pii/S2352467723002126},
author = {Marta Rodrigues and Tiago Soares and Hugo Morais}
}

@article{nawaz2020tso,
author = {Nawaz, Aamir and Wang, Hongtao and Wu, Qiuwei and Kumar Ochani, Manesh},
title = {TSO and DSO with large-scale distributed energy resources: A security constrained unit commitment coordinated solution},
journal = {International Transactions on Electrical Energy Systems},
volume = {30},
number = {3},
pages = {e12233},
doi = {https://doi.org/10.1002/2050-7038.12233},
url = {https://onlinelibrary.wiley.com/doi/abs/10.1002/2050-7038.12233},
eprint = {https://onlinelibrary.wiley.com/doi/pdf/10.1002/2050-7038.12233},
note = {e12233 ITEES-19-0196.R3},
year = {2020}
}

@article{HAJATI2024110840,
title = {Maximizing social welfare in local flexibility markets by integrating the value of flexibility loss (VOFL)},
journal = {Electric Power Systems Research},
volume = {235},
pages = {110840},
year = {2024},
issn = {0378-7796},
doi = {https://doi.org/10.1016/j.epsr.2024.110840},
url = {https://www.sciencedirect.com/science/article/pii/S0378779624007260},
author = {Maryam Hajati and Mohamad Kazem Sheikh-El-Eslami and Hamed Delkhosh}
}

@ARTICLE{9006872,
  author={Evangelopoulos, Vasileios A. and Avramidis, Iason I. and Georgilakis, Pavlos S.},
  journal={IEEE Access}, 
  title={Flexibility Services Management Under Uncertainties for Power Distribution Systems: Stochastic Scheduling and Predictive Real-Time Dispatch}, 
  year={2020},
  volume={8},
  number={},
  pages={38855-38871},
  doi={10.1109/ACCESS.2020.2975663}}

@article{KALANTARNEYESTANAKI2024110747,
title = {Grid-cognizant TSO and DSO coordination framework for active and reactive power flexibility exchange: The Swiss case study},
journal = {Electric Power Systems Research},
volume = {235},
pages = {110747},
year = {2024},
issn = {0378-7796},
doi = {https://doi.org/10.1016/j.epsr.2024.110747},
url = {https://www.sciencedirect.com/science/article/pii/S0378779624006333},
author = {Mohsen Kalantar-Neyestanaki and Rachid Cherkaoui}
}

@ARTICLE{9939101,
  author={Bakhtiari, Hamed and Hesamzadeh, Mohammad Reza and Bunn, Derek W.},
  journal={IEEE Transactions on Power Systems}, 
  title={TSO-DSO Operational Coordination Using a Look-Ahead Multi-Interval Framework}, 
  year={2023},
  volume={38},
  number={5},
  pages={4221-4239},
  doi={10.1109/TPWRS.2022.3219581}}

@article{LIU202327,
title = {To exploit the flexibility of TSO–DSO interaction: A coordinated transmission robust planning and distribution stochastic reinforcement solution},
journal = {Energy Reports},
volume = {9},
pages = {27-36},
year = {2023},
note = {2022 9th International Conference on Power and Energy Systems Engineering},
issn = {2352-4847},
doi = {https://doi.org/10.1016/j.egyr.2022.10.368},
url = {https://www.sciencedirect.com/science/article/pii/S2352484722023034},
author = {Jia Liu and Libo Zhang and Kaicheng Liu and Zao Tang and Peter Pingliang Zeng and Yalou Li}
}

@article{MANSOURI2023121062,
title = {An interval-based nested optimization framework for deriving flexibility from smart buildings and electric vehicle fleets in the TSO-DSO coordination},
journal = {Applied Energy},
volume = {341},
pages = {121062},
year = {2023},
issn = {0306-2619},
doi = {https://doi.org/10.1016/j.apenergy.2023.121062},
url = {https://www.sciencedirect.com/science/article/pii/S0306261923004269},
author = {Seyed Amir Mansouri and Emad Nematbakhsh and Ahmad Rezaee Jordehi and Mousa Marzband and Marcos Tostado-Véliz and Francisco Jurado}
}

@article{CHEN2022118319,
title = {Optimal participation of ADN in energy and reserve markets considering TSO-DSO interface and DERs uncertainties},
journal = {Applied Energy},
volume = {308},
pages = {118319},
year = {2022},
issn = {0306-2619},
doi = {https://doi.org/10.1016/j.apenergy.2021.118319},
url = {https://www.sciencedirect.com/science/article/pii/S0306261921015749},
author = {Houhe Chen and Di Wang and Rufeng Zhang and Tao Jiang and Xue Li}
}

@article{SOTO2021116268,
title = {Peer-to-peer energy trading: A review of the literature},
journal = {Applied Energy},
volume = {283},
pages = {116268},
year = {2021},
issn = {0306-2619},
doi = {https://doi.org/10.1016/j.apenergy.2020.116268},
url = {https://www.sciencedirect.com/science/article/pii/S0306261920316585},
author = {Esteban A. Soto and Lisa B. Bosman and Ebisa Wollega and Walter D. Leon-Salas}
}

@article{MARQUES2023101055,
title = {P2P flexibility markets models to support the coordination between the transmission system operators and distribution system operators},
journal = {Sustainable Energy, Grids and Networks},
volume = {34},
pages = {101055},
year = {2023},
issn = {2352-4677},
doi = {https://doi.org/10.1016/j.segan.2023.101055},
url = {https://www.sciencedirect.com/science/article/pii/S2352467723000632},
author = {João Marques and Tiago Soares and Hugo Morais}
}

@article{LI2024122328,
title = {An electricity and carbon trading mechanism integrated with TSO-DSO-prosumer coordination},
journal = {Applied Energy},
volume = {356},
pages = {122328},
year = {2024},
issn = {0306-2619},
doi = {https://doi.org/10.1016/j.apenergy.2023.122328},
url = {https://www.sciencedirect.com/science/article/pii/S0306261923016926},
author = {Junkai Li and Shaoyun Ge and Hong Liu and Qi Yu and Shida Zhang and Chengshan Wang and Chenghong Gu}
}

@ARTICLE{9911669,
  author={Ullah, Md Habib and Park, Jae-Do},
  journal={IEEE Transactions on Power Systems}, 
  title={Transactive Energy Market Operation Through Coordinated TSO-DSOs-DERs Interactions}, 
  year={2023},
  volume={38},
  number={2},
  pages={1978-1990},
  doi={10.1109/TPWRS.2022.3212065}}

@misc{zemkohoo2020bilevel,
  title={Bilevel optimization advances and next challenges},
  author={Zemkohoo, Alain and Dempe, Stephan},
  year={2020},
  publisher={Springer Geneva, Switzerland}
}

@ARTICLE{4956966,
  author={Stott, Brian and Jardim, Jorge and Alsac, Ongun},
  journal={IEEE Transactions on Power Systems}, 
  title={DC Power Flow Revisited}, 
  year={2009},
  volume={24},
  number={3},
  pages={1290-1300},
  doi={10.1109/TPWRS.2009.2021235}}

@article{cain2012history,
  title={History of optimal power flow and formulations},
  author={Cain, Mary B and O’neill, Richard P and Castillo, Anya and others},
  journal={Federal Energy Regulatory Commission},
  volume={1},
  pages={1--36},
  year={2012},
  publisher={Citeseer}
}

@book{wood2013power,
  title={Power generation, operation, and control},
  author={Wood, Allen J and Wollenberg, Bruce F and Shebl{\'e}, Gerald B},
  year={2013},
  publisher={John wiley \& sons}
}

@ARTICLE{1664986,
  author={Jabr, R.A.},
  journal={IEEE Transactions on Power Systems}, 
  title={Radial distribution load flow using conic programming}, 
  year={2006},
  volume={21},
  number={3},
  pages={1458-1459},
  doi={10.1109/TPWRS.2006.879234}}

@ARTICLE{7990560,
  author={Molzahn, Daniel K. and Dörfler, Florian and Sandberg, Henrik and Low, Steven H. and Chakrabarti, Sambuddha and Baldick, Ross and Lavaei, Javad},
  journal={IEEE Transactions on Smart Grid}, 
  title={A Survey of Distributed Optimization and Control Algorithms for Electric Power Systems}, 
  year={2017},
  volume={8},
  number={6},
  pages={2941-2962},
  doi={10.1109/TSG.2017.2720471}}

@article{GIVISIEZ2020106659,
title = {A Review on TSO-DSO Coordination Models and Solution Techniques},
journal = {Electric Power Systems Research},
volume = {189},
pages = {106659},
year = {2020},
issn = {0378-7796},
doi = {https://doi.org/10.1016/j.epsr.2020.106659},
url = {https://www.sciencedirect.com/science/article/pii/S0378779620304624},
author = {Arthur Gonçalves Givisiez and Kyriacos Petrou and Luis F. Ochoa}
}

@article{ZHANG2024123073,
title = {Establishing a hierarchical local market structure using multi-cut Benders decomposition},
journal = {Applied Energy},
volume = {363},
pages = {123073},
year = {2024},
issn = {0306-2619},
doi = {https://doi.org/10.1016/j.apenergy.2024.123073},
url = {https://www.sciencedirect.com/science/article/pii/S0306261924004562},
author = {Haoyang Zhang and Sen Zhan and Koen Kok and Nikolaos G. Paterakis}
}

@article{YING2024123803,
title = {Decentralized energy management of a hybrid building cluster via peer-to-peer transactive energy trading},
journal = {Applied Energy},
volume = {372},
pages = {123803},
year = {2024},
issn = {0306-2619},
doi = {https://doi.org/10.1016/j.apenergy.2024.123803},
url = {https://www.sciencedirect.com/science/article/pii/S0306261924011863},
author = {Chenhao Ying and Yunyang Zou and Yan Xu}
}

@article{HUANG2022108179,
title = {Bilateral energy-trading model with hierarchical personalized pricing in a prosumer community},
journal = {International Journal of Electrical Power \& Energy Systems},
volume = {141},
pages = {108179},
year = {2022},
issn = {0142-0615},
doi = {https://doi.org/10.1016/j.ijepes.2022.108179},
url = {https://www.sciencedirect.com/science/article/pii/S0142061522002137},
author = {Ting Huang and Yi Sun and Mengting Jiao and Zhuang Liu and Jianhong Hao}
}

@incollection{shahidehpour2003ieee,
  title={IEEE-30 Bus System Data},
  author={Shahidehpour, Mohammad and Wang, Yaoyu and others},
  booktitle={Communication and control in electric power systems: applications of parallel and distributed processing},
  pages={493--495},
  year={2003}
}

@ARTICLE{9258930,
  author={Dolatabadi, Sarineh Hacopian and Ghorbanian, Maedeh and Siano, Pierluigi and Hatziargyriou, Nikos D.},
  journal={IEEE Transactions on Power Systems}, 
  title={An Enhanced IEEE 33 Bus Benchmark Test System for Distribution System Studies}, 
  year={2021},
  volume={36},
  number={3},
  pages={2565-2572},
  doi={10.1109/TPWRS.2020.3038030}}

@ARTICLE{9328796,
  author={Sun, Xianzhuo and Qiu, Jing},
  journal={IEEE Transactions on Smart Grid}, 
  title={Two-Stage Volt/Var Control in Active Distribution Networks With Multi-Agent Deep Reinforcement Learning Method}, 
  year={2021},
  volume={12},
  number={4},
  pages={2903-2912},
  doi={10.1109/TSG.2021.3052998}}

@article{zemkoho2021theoretical,
  title={Theoretical and numerical comparison of the Karush--Kuhn--Tucker and value function reformulations in bilevel optimization},
  author={Zemkoho, Alain B and Zhou, Shenglong},
  journal={Computational Optimization and Applications},
  volume={78},
  number={2},
  pages={625--674},
  year={2021},
  publisher={Springer}
}

@article{dempe2015bilevel,
  title={Bilevel programming problems},
  author={Dempe, Stephan and Kalashnikov, Vyacheslav and P{\'e}rez-Vald{\'e}s, Gerardo A and Kalashnykova, Nataliya},
  journal={Energy Systems. Springer, Berlin},
  volume={10},
  number={978-3},
  pages={53--56},
  year={2015},
  publisher={Springer}
}

@ARTICLE{1216140,
  author={Bakirtzis, A.G. and Biskas, P.N.},
  journal={IEEE Transactions on Power Systems}, 
  title={A decentralized solution to the DC-OPF of interconnected power systems}, 
  year={2003},
  volume={18},
  number={3},
  pages={1007-1013},
  doi={10.1109/TPWRS.2003.814853}}
\end{document}